\documentclass[reqno,final]{siamltex}
\usepackage[shortlabels]{enumitem}
\usepackage{geometry}
\usepackage{layouts}
\usepackage{booktabs}
\usepackage{stmaryrd}
\usepackage{harpoon}
\usepackage{xcolor}
\usepackage{graphicx}
\usepackage{subfigure}
\usepackage{cite}
\usepackage{diagbox}
\usepackage{mathrsfs}
\usepackage{amsmath,amssymb,amscd,amsxtra,amsfonts}
\usepackage{algorithm,algorithmic}
\usepackage{comment}
\usepackage{bm}
\usepackage{multirow}
\usepackage{ulem}

\usepackage[colorlinks,
linkcolor=blue,
anchorcolor=red,
citecolor=green]{hyperref}
\usepackage{tikz}

\providecommand{\Div}{\operatorname{div}}          
\providecommand{\curl}{\operatorname{{\bf curl}}}  

\providecommand*{\Dist}[2]{\operatorname{dist}({#1};{#2})}   
\providecommand*{\Dist}[2]{\Dist{#1}{#2}}
\providecommand*{\Span}[1]{\operatorname{Span}\left\{{#1}\right\}}     

\newcommand{\Ve}{{\mathbf{e}}}
\newcommand{\Vf}{{\mathbf{f}}}
\newcommand{\Vg}{{\mathbf{g}}}

\newcommand{\Vr}{{\mathbf{r}}}

\newcommand{\Vu}{{\mathbf{u}}}
\newcommand{\Vv}{{\mathbf{v}}}
\newcommand{\Vw}{{\mathbf{w}}}

\newcommand{\Vz}{{\mathbf{z}}}

\newcommand{\Bn}{{\boldsymbol{n}}}

\newcommand{\Bx}{{\boldsymbol{x}}}

\newcommand{\BF}{{\boldsymbol{F}}}

\newcommand{\BL}{{\boldsymbol{L}}}

\newcommand{\xibf}{\boldsymbol{\xi}}

\newcommand{\sigmabf}{\boldsymbol{\sigma}}

\newcommand{\varphibf}{\boldsymbol{\varphi}}

\newcommand{\Ca}{{\cal A}}

\newcommand{\Cf}{{\cal F}}

\newcommand{\Cm}{{\cal M}}
\newcommand{\Cn}{{\cal N}}
\newcommand{\Co}{{\cal O}}
\newcommand{\Cp}{{\cal P}}

\newcommand{\Cr}{{\cal R}}
\newcommand{\Cs}{{\cal S}}
\newcommand{\Ct}{{\cal T}}

\newcommand{\bbA}{\mathbb{A}}

\newcommand{\bbC}{\mathbb{C}}
\newcommand{\bbD}{\mathbb{D}}
\newcommand{\bbE}{\mathbb{E}}

\newcommand{\bbG}{\mathbb{G}}

\newcommand{\bbI}{\mathbb{I}}

\newcommand{\bbM}{\mathbb{M}}

\newcommand{\bbP}{\mathbb{P}}

\newcommand{\bbR}{\mathbb{R}}

\newcommand*{\DP}[2]{\left<{#1},{#2}\right>} 

\newcommand*{\N}[1]{\|{#1}\|}     
\newcommand*{\SN}[1]{|{#1}|}      
\newcommand*{\SP}[2]{({#1},{#2})} 

\newcommand*{\Lp}[2][\defaultdomain]{L^{#2}({#1})}
\newcommand*{\Lpv}[2][\defaultdomain]{\BL^{#2}({#1})}
\newcommand*{\NLp}[3][\defaultdomain]{\N{#2}_{\Lp[#1]{#3}}}

\newcommand*{\Ltwo}[1][\defaultdomain]{\Lp[#1]{2}}
\newcommand*{\Ltwov}[1][\defaultdomain]{\Lpv[#1]{2}}
\newcommand*{\NLtwo}[2][\defaultdomain]{\NLp[#1]{#2}{2}}

\newcommand*{\Hm}[2][\defaultdomain]{H^{#2}({#1})}

\newcommand*{\Hdiv}[1][\defaultdomain]{\boldsymbol{H}(\Div,{#1})}

\newcommand*{\Hcurl}[1][\defaultdomain]{\boldsymbol{H}(\curl,{#1})}

\newcommand*{\jump}[2][]{ \llbracket{#2}\rrbracket_{#1}}

\newcommand*{\avg}[2][]{\{\hskip -3.5pt\{{#2}\}\hskip -3.5pt\}_{#1}}

\newcommand{\be}{\begin{eqnarray}}
\newcommand{\ee}{\end{eqnarray}}

\newcommand{\ben}{\begin{eqnarray*}}
\newcommand{\een}{\end{eqnarray*}}

\newtheorem{remark}[theorem]{\sc Remark}

\def\address#1{\expandafter\def\expandafter\@aabuffer\expandafter
	{\@aabuffer{\affiliationfont{#1}}\relax\par
	\vspace*{13pt}}}

\title{\lowercase{$p$}-multigrid method for the discontinuous Galerkin discretization of elliptic problems}

\author{
  Nuo Lei \thanks{
    Hua Loo-Keng Center for Mathematical Sciences,
    Academy of Mathematics and Systems Science, 
    Chinese Academy of Sciences, 
    Beijing, 100190, China. 
    The first author was supported in part by 
    the Strategic Priority Research Program of the Chinese Academy of Sciences XDB0640000 and XDB0640100,
    National Natural Science Foundation of China 12288201. 
    (nuo\_lei@lsec.cc.ac.cn)
  }
  \and Donghang Zhang \thanks{
    Academy of Mathematics and Systems Science, 
    Chinese Academy of Sciences, 
    Beijing, 100190, China. 
    The second author was supported in part by 
    the Strategic Priority Research Program of the Chinese Academy of Sciences XDB0640000 and XDB0640300,
    the Postdoctoral Fellowship Program (Grade B) of China Postdoctoral Science Foundation GZB20230816, 
    and National Natural Science Foundation of China 12401487.
    (zdh@lsec.cc.ac.cn)
  }
  \and Weiying Zheng \thanks{
    LSEC, Institute of Computational Mathematics and Scientific/Engineering Computing, 
    Academy of Mathematics and Systems Science, 
    Chinese Academy of Sciences; 
    School of Mathematical Science,
    University of Chinese Academy of Sciences, 
    Beijing, 100190, China. 
    The third author was supported in part by 
    the Strategic Priority Research Program of the Chinese Academy of Sciences XDB0640000 and XDB0640100,
    National Key R \& D Program of China 2019YFA0709600 and 2019YFA0709602. 
    (zwy@lsec.cc.ac.cn)
  }
}

\begin{document}
\maketitle

\begin{abstract}
In this paper, we propose a $W$-cycle $p$-multigrid method 
for solving the $p$-version symmetric interior penalty discontinuous Galerkin (SIPDG) 
discretization of elliptic problems. 
This SIPDG discretization employs hierarchical Legendre polynomial basis functions. 
Inspired by the uniform convergence theory of the $W$-cycle $hp$-multigrid method in 
[P. F. Antonietti, et al., SIAM J. Numer. Anal., 53 (2015)], 
we provide a rigorous convergence analysis 
for the proposed $p$-multigrid method, 
considering both inherited and non-inherited bilinear forms of SIPDG discretization. 
Our theoretical results show significant improvement over 
[P. F. Antonietti, et al., SIAM J. Numer. Anal., 53 (2015)], 
reducing the required number of smoothing steps from $\Co(p^2)$ to $\Co(p)$, 
where $p$ is the polynomial degree of the discrete broken polynomial space. 
Moreover, the convergence rate remains independent of the mesh size.
Several numerical experiments are presented to verify our theoretical findings.
Finally,  we numerically verify the effectiveness of the $p$-multigrid method for unfitted finite element discretization in solving elliptic interface problems on general $C^{2} $-smooth interfaces.

\end{abstract}



\section{Introduction}
\label{sec1}
The Poisson equation is a fundamental partial differential equation (PDE) 
that arises in numerous scientific and engineering disciplines. 
Efficient solvers are relatively mature for low-order approximations, 
such as geometric multigrid (GMG) methods \cite{Xu1992} and algebraic multigrid (AMG) methods \cite{Xu2017_review}. 
However, high-order methods are often required to capture detailed features of the solution, 
especially in complex geometries or interface problems,
and they may lead to large linear systems that are computationally expensive to solve.

Discontinuous Galerkin (DG) methods have emerged as a powerful and versatile class of numerical techniques. 
They use local basis functions on each element and allow discontinuities between elements, 
making them particularly suitable for handling complex geometries and interface problems.
Additionally, DG methods can be easily extended to high-order schemes 
by increasing the polynomial degree of the basis functions within each element.
The interior penalty DG (IPDG) method is a specific variant of DG method which incorporates penalty terms to enforce continuity between elements \cite{Arnold1982, Arnold2002}. 
The development of IPDG methods is driven by the need to improve stability 
and convergence properties while maintaining the inherent flexibility and high-order accuracy of the DG methods.

Despite the advantages mentioned above, 
the IPDG method poses significient challenges in solving the resulting discrete problem.
High-order discretizations lead to large and ill-conditioned linear systems, 
making traditional direct solvers impractical due to their 
inefficiency and memory demands.
Iterative solvers are more suitable for these large systems, 
they require effective preconditioners or advanced techniques 
like multigrid methods or domain decomposition methods (DDM) to ensure convergence.

Developing such efficient solvers is crucial for algebraic systems arising from high-order IPDG methods. 
A preconditioner effectively transforms the stiffness matrix into one with a better condition number, 
clustering its spectrum and thereby accelerating the convergence of iterative methods. 
This not only enhances the overall computational efficiency 
but also ensures that the solver remains robust to mesh sizes and problem configurations.
In \cite{Bastian2012}, the authors proposed an AMG method based on subspace corrections,
which is robust and only mildly affected by coefficient jumps and mesh sizes.
Another AMG method is proposed in \cite{Antonietti2020}
and shows uniform convergence with respect to mesh sizes and polynomial degrees. 
In \cite{Hartmann2009}, the authors 
derived and analyzed a block diagonal preconditioner, 
that is robust to mesh sizes for IPDG discretization. 
Based on space decompositions or auxiliary space methods, 
the authors in \cite{Dios2009,Dios2014,Dios2017} proposed a family of 
robust preconditioners for elliptic problems with jump coefficients.
A robust additive Schwarz preconditioner was developed for high-order discretizations \cite{Antonietti2017_JSC}, 
and the preconditioner's spectral radius 
is shown to be invariant with respect to key discretization parameters.
In \cite{Thiele2017}, the authors proposed an inexact hierarchical scale separation (HSS) framework that 
dynamically adapts coarse-solver tolerances based on fine-scale residuals,
and achieves significant speedups for the Cahn-Hilliard equation.

Multigrid methods are efficient numerical methods for solving linear systems, 
arising from finite element or finite difference discretization of elliptic problems \cite{Xu1992}. 
They can be categorized into GMG and AMG methods. 
GMG methods use a predefined sequence of grids based on the problem's geometry, 
whereas AMG methods construct the hierarchy of grids algebraically 
from the stiffness matrix of the linear system, 
making them more flexible for complex or unstructured meshes.
Some geometric multigrid methods have been applied to the 
preconditioning of DG systems \cite{Brenner2005, Lazarov2006,Lazarov2007,Xu2024}.
Recently, there have been some GMG works on unstructured polytopic meshes.
The framework in \cite{Antonietti2019} establishes uniform $hp$-convergence 
on agglomerated meshes via non-nested $L^2$-projections and adaptive Schwarz smoothing.
Using physics-preserving Dirichlet-to-Neumann intergrid transfers in \cite{Wildey2019}, 
the proposed methods avoid parameter upscaling and demonstrate robust scalability for HDG/IPDG methods.

The $p$-multigrid ($p$-MG) method is a variant of the standard multigrid method 
that exploits polynomial-degree hierarchy instead of mesh hierarchy.
In contrast, the $p$-MG method maintains a fixed mesh,
thus being suitable for complex geometries.
The $p$-MG approach is especially beneficial for spectral method \cite{Ronquist1987,Maday1988} 
and high-order discretizations, 
where increasing polynomial degrees can significantly improve the accuracy without refining meshes.
In \cite{Helenbrook2003,Helenbrook2006},
the authors investigated the $p$-MG method for DG discretization of the Poisson equation.
By combining $p$-MG methods with AMG techniques, 
the authors proposed an efficient solver for pore scale flow in \cite{Thiele2022}. 
An improved $p$-MG method was proposed by Gunatilake, 
that uses a global smoothing strategy via space decomposition in a multiplicative Schwarz framework \cite{Gunatilake2022}. 
The method demonstrates better convergence rates than existing $p$-MG methods.
We also refer to $p$-MG applications for the Euler equations \cite{Luo2006},
the Helmholtz equations \cite{Haupt2013}, the Stokes equations, 
and linear elasticity equations \cite{Hong2016}.

Although $p$-MG methods of high-order DG discretizations have been widely used 
and have been shown to be highly efficient through numerical experiments, 
rigorous theoretical analysis is quite limited.
In \cite{Antonietti2015,Antonietti2017}, an $hp$-multigrid method was proposed for elliptic problems
with rigorous convergence analysis.
To the best of our knowledge, this is the first theoretical result for $p$- and $hp$-DG methods.
Our work is inspired by \cite{Antonietti2015}, 
and improves these results with a better convergence rate.

This paper is our first step toward studying $p$-MG for DG-discretized problems.
We investigate a $p$-MG method for solving the SIPDG discretization of elliptic problems. 
The contributions of the paper are listed as follows.
\begin{itemize}
\item 
We adopt the fourth-kind Chebyshev polynomial iterative method \cite{Lottes2023} as smoother of the $p$-MG method,
and prove pre- and post-smoothing properties of the polynomial smoother.
Numerical experiments validate our theoretical analysis.
\item
We present a rigorous convergence analysis of the $W$-cycle $p$-MG method, 
which uses both inherited and non-inherited bilinear forms.
The convergence rate remains independent of the mesh size, and requires $\Co(p)$
smoothing steps. 
\item 
We provide extensive numerical results in 1D, 2D and 3D
and compare our $p$-MG method with $p$-MG methods in the literature.
\end{itemize}
While implemented on Cartesian meshes 
(providing a controlled testbed for analyzing approximation/smoothing properties), 
the core techniques can generalize to arbitrary mesh topologies.

The rest of the paper is organized as follows.
In Section \ref{sec2}, 
we introduce high-order SIPDG discretization of Poisson equation. 
In Section \ref{sec3}, we present the $p$-MG algorithm 
using non-inherited bilinear forms, 
the transfer operators 
and the polynomial smoother.
The convergence analysis of two-level 
and $W$-cycle $p$-MG algorithms are given in Section \ref{sec4}, 
with approximation and smoothing properties.
Section \ref{sec5} details the inherited $p$-MG algorithm and its analysis.
Various numerical experiments are given in Section \ref{sec6} and Section \ref{sec7}
to validate the excellent performance of our methods.
Finally, Section \ref{sec8} provides a summary of the method.

\section{Model problem and IPDG discretization}
\label{sec2}
Suppose $D \subset \mathbf{R}^{d} $ is an open and bounded domain. 
Let $\Ltwo[D]$ be the space of square-integrable functions. 
The inner product and norm on $\Ltwo[D]$ are defined by 
 \[
   \SP{u}{v}_{D} := \int_{D} u v, \quad 
   \NLtwo[D]{u} := \sqrt{\SP{u}{u}_{D} }.
\] 
For $(d-1)$-dimensional manifold $F \subset \partial D $, the inner product and norm on $\Ltwo[F]$ are denoted by
\[
   \DP{u}{v}_{F} := \int_{F} u v, \quad
   \NLtwo[F]{u} := \sqrt{\DP{u}{u}_{F} }.
\] 
For any integer $s > 0$, denote 
$\Hm[D]{s}:=\{ v \in \Ltwo[D] : \N{v}_{H^s(D)} < +\infty \}$,
with $\N{v}_{H^s(D)} := ( \sum_{\N{\sigmabf}_{\infty} \leq s } \int_D |\partial^{\sigmabf} v(\Bx) |^2 d \Bx )^{1/2}$
where $\sigmabf$ is a non-negative $d$-tuple, and $\|\sigmabf\|_{\infty} := \max\{\sigma_{i}:i=1,\cdots,d\}  $. 
Moreover, vector-valued quantities will be denoted by boldface notations, such as $\Ltwov[D] = \left(\Ltwo[D]\right)^{d}, 
\sigmabf= (\sigma_1, \cdots, \sigma_d ) ,$ etc.

In this section, we introduce the model problem and its discretization by 
$p$-version IPDG formulations.
We describe the IPDG formulations using hierarchical bases, 
which will play a pivotal role in the forthcoming $p$-MG analysis.
In this paper, we suppose $\Omega \subset \mathbf{R}^{d} $ is rectangular domain, 
and consider the following model problem, 
\begin{equation}\label{cont}
-\Delta u = f \quad \text{in}\;\; \Omega, \quad
u=0\quad \text{on}\;\;  \partial \Omega,
\end{equation}
where $f \in \Ltwo$ is a given function.

\subsection{The IPDG discretization}
\label{sec2.1}
Suppose $\mathcal{T}_h$ is a Cartesian uniform mesh of $\Omega$
consisting of squares or cuboid elements with edge size $h$,
and the sides of an element $T$ are parallel to the coordinate axes. 
The total number of elements in $\Ct_h$ is denoted by $N_h$. 
Set $\Cf_h = \Cf_h^i \bigcup \Cf_h^b$, 
where $\Cf_h^i$ and $\Cf_h^b$ denote the sets of all interior faces and boundary faces, respectively.

The Legendre polynomials on the reference element $\widehat{T} = [-1,1]^{d} $ are denoted by
\[
{L}_{\sigmabf} (\xibf) = \prod_{i=1} ^{d} l_{\sigma_{i}} (\xi_i), \quad
\sigmabf \in \mathbf{N}^{d} , \xibf \in \widehat{T},
\]
where $\sigmabf = (\sigma_1, \cdots, \sigma_d)^\top$ and $\xibf = (\xi_1,\cdots, \xi_d)^\top$ are 
$d$-dimensional tuples with 
$\sigma_i \in \mathbf{N}^+, \xi_i \in [-1, 1] $,
and $l_{\sigma_i} (\xi_i)$ is the Legendre polynomial of degree $\sigma_i$
for any $i = 1,\cdots, d$.
For each element $T\in \Ct_h$, let $\Bx=\BF_T(\xibf)$ denote the affine transformation from $\widehat{T}$ to $T$. 
The Legendre polynomials on $T$ are given by
\begin{align}\label{basis}
\Psi^{\sigmabf}_T(\Bx) :=  L_{\sigmabf}(\BF_T^{-1}(\Bx)) = \prod_{i=1}^d \psi_T^{\sigma_i} (x_i), \quad \text{with} \quad
\psi^{\sigma_i}_T(x_i) := l_{\sigma_i}(\xi_i).
\end{align}
For $k = 1, 2, \dots, p$, the space of polynomials with degree no more than $k$ for each variable is denoted 
\[
  Q_{k} (T) : = \Span{\Psi^{\sigmabf}_T : \|\sigmabf\|_\infty \le k \, \forall \, \sigmabf \in \mathbf{N}^{d} },
  \quad 
  N_k := \text{dim}  Q_{k} (T) = (k+1)^d . 
\]
The broken polynomial space is defined as
\begin{align}\label{DGspace}
U_k(\Ct_h) :=
\{v\in\Ltwo: v|_T\in Q_k(T) \, \forall \, T \in \Ct_{h} \}.
\end{align}
We obtain the space hierarchy
\begin{align*}
U_{1} (\Ct_{h}) \subset U_{2} (\Ct_{h} ) \subset \cdots \subset U_{p} (\Ct_{h}). 
\end{align*} 
To each $F\in \Cf_h^i$ with $F = T_1\cap T_2$, $T_1\neq T_2$, 
the average and jump of a scalar-valued function $u$ are defined on $F$ by,
\ben
\avg[F]{u} := \frac12u|_{T_1} + \frac12u|_{T_2}, \quad
\jump[F]{u} := u|_{T_1} - u|_{T_2}.
\een 
Moreover, we keep the orientation of unit normal $\Bn_{F} $ is arbitrary but fixed from $T_{1} $ to $T_{2} $. 
If $F\in \Cf_h^b$ is a boundary face, 
$\avg[F]{u} = \jump[F]{u} := u$, 
and $\Bn_{F}$ is defined as the unit outward normal to $\partial \Omega$.
Without confusions, we omit the subscript $F$ 
and simply write $\avg{u}$, $\jump{u}$, and $\Bn$. 
For each $k = 1, \dots, p$, we define the symmetric bilinear form 
$\Ca_{k} (\cdot, \cdot): U_{k}(\Ct_{h}) \times U_{k} (\Ct_{h}) \to \mathbf{R}$ 
\begin{align}\label{IPDG-a}
\begin{split}
\Ca_{k} (u, v) :=& 
\sum_{T \in \Ct_{h} } \SP{\nabla u}{\nabla v}_{T}  
- \sum_{F \in \Cf_{h} }\DP{\avg{ \nabla u }\cdot\Bn }{\jump{v }}_{F} \\
&- \sum_{F\in\Cf_{h} } \DP{\avg{\nabla v }\cdot \Bn }{\jump{u} }_{F}
+ \sum_{F \in \Cf_{h} } \alpha \DP{\jump{u}}{\jump{v}}_{F}, 
\end{split}
\end{align}
where $\alpha := {\alpha_0 k^{2} }{h}^{-1} $ and $\alpha_0 > 0$ is a constant large enough \cite{Epshteyn2007}.
The $p$-version IPDG formulation is to find $u_k  \in U_k(\Ct_h)$ such that
\begin{align}\label{weak-h}
\Ca_{k} (u_{k} , v_{k} ) = (f, v_{k} )_{\Omega}  \quad \forall \, v_k \in U_k(\Ct_h).
\end{align}
For any $u_k \in U_{k} (\Ct_{h})$ and $1\leq k\leq p$, the expansion holds
\begin{align}\label{basisexp}
u_k =  \sum_{T \in \Ct_{h} } \sum_{\|\sigmabf\|_{\infty} \le k } u_{T}^{\sigmabf} \Psi_{T}^{\sigmabf} = 
\sum_{T \in \Ct_{h} } \sum_{\|\sigmabf\|_{\infty} \le k } u_{T}^{\sigmabf} \prod_{i=1}^{d} \psi_{T} ^{\sigma_i} .
\end{align}
Throughout the paper, we use $\Vu_{k} \in \mathbf{R}^{N_{h} N_{k}}$
to denote the vector of coefficients $\{u_{T} ^{\sigmabf} \in \mathbf{R} : \N{\sigmabf}_{\infty} \le k, T \in \Ct_{h}  \}$. 
Clearly, \eqref{weak-h} is equivalent to the algebraic equations, 
\begin{equation}\label{LinearSystem}
	\bbA_k \Vu_k = \Vf_k,
\end{equation}
where $\bbA_k \in \mathbf{R}^{N_{h} N_{k} \times N_{h} N_{k}  } $ 
and $\Vf_{k} \in \mathbf{R}^{N_{h} N_{k} } $ are the stiffness matrix 
and the load vector, respectively.
Before describing our $p$-MG method, 
we recall some technical results that will be used in the forthcoming analysis.

\subsection{Technical results}
\label{sec2.2}
Throughout this paper, $x\lesssim y$ indicates $x \leq Cy$ for a constant $C > 0$, 
independent of the mesh size and polynomial degrees. 
In the meantime, $x\approx y$ denotes that there exist constants $c, C > 0$, 
independent of the mesh size and polynomial degrees, such that $cy \leq x \leq Cy$.
We endow the broken polynomial space $U_{k} (\Ct_{h})$ with the mesh-dependent norm $\N{\cdot}_{1,k,h}$ defined as
\[
\N{v}_{1,k,h}^{2}  := 
\sum_{T\in\Ct_{h}} \NLtwo[T]{\nabla v}^{2} + 
\frac{ k^{2}  }{ h } \sum_{F\in\Cf_{h}} \NLtwo[F]{\jump{v}}^{2} . 
\] 
The following lemma ensures the bilinear form $\Ca_{k} (\cdot, \cdot)$ is continuous and coercive under the mesh-dependent norm $\N{\cdot}_{1,k,h} $, see \cite[Proposition 3.1]{Perugia2002}.
\begin{lemma}\label{lem2.1}
Suppose $\alpha_0$ is large enough. 
For any level $k = 1,\dots, p$, it holds that
\begin{align}
\Ca_{k} (v_k,v_k) \gtrsim \N{v_k}^2_{1,k,h},\quad
\SN{\Ca_{k} (v_k,w_k)} \lesssim \N{v_k}_{1,k,h}\N{w_k}_{1,k,h} \quad
\forall \, v_{k}, w_k \in U_{k} (\Ct_{h}) .
\end{align}
\end{lemma}

The broken Sobolev space is denoted by $\Hm[\Ct_{h}]{s}, s \ge 1$ which consists of all element-wise $H^{s}$-functions.
Moreover, we have the following optimal error estimates 
with homogeneous Dirichlet boundary conditions, cf. \cite[Theorem 3.3]{Stamm10}, and \cite[Remark 2.3]{Antonietti2015}.
\begin{lemma} 
Suppose $u_k \in U_{k} (\Ct_{h})$ is the solution to the discrete problem, 
\begin{equation}\label{stiffness_Ak}
  \Ca_{k} (u_{k} ,v_{k} ) = \SP{f}{v_{k} }_{\Omega} \quad \forall \, v_k\in U_k(\Ct_h) .   
\end{equation}
The following error estimates hold, for $u\in H^{s+1}(\Ct_{h} )$ and $s \geq 1$,
\begin{align}
& \N{u-u_k}_{1,k,h} \lesssim \frac{h^{\min(k,s)}}{k^{s}} \N{u}_{H^{s+1}(\Ct_h)} ,\label{error_estimate_1}\\
& \N{u-u_k}_{L^2(\Omega)} \lesssim \frac{h^{\min(k,s)+1}}{k^{s+1}} \N{u}_{H^{s+1}(\Ct_h)} .\label{error_estimate_2}
\end{align}
\end{lemma}


\section{$p$-MG method with non-inherited bilinear forms}
\label{sec3}
In this section, we will describe the $W$-cycle algorithm of the $p$-MG method with non-inherited bilinear forms.
We first introduce the restriction and prolongation operators between broken polynomial spaces,
then introduce the polynomial smoother and the $W$-cycle algorithm.

For any $v_k, w_k \in U_k(\Ct_h)$, 
their vectors of unknowns are denoted by $\Vv_{k} $ and $\Vw_{k} $, respectively, 
according to \eqref{basisexp}.
The matrix representations of the non-inherited bilinear forms $\{\Ca_{k}(\cdot,\cdot)\}_{k=1} ^{p} $ are defined as
\begin{align*}
\Vv_{k} ^\top\bbA_k\Vw_{k}:=
\Ca_{k} (v_{k} , w_{k}) \quad \forall \, v_{k} , w_{k}  \in U_{k} (\Ct_{h}).
\end{align*}
By Lemma \ref{lem2.1}, matrix $\bbA_{k} $ is symmetric and positive definite. 
Then we can define the $\bbA_{k}$-dependent inner product
and the $s$-norm ($s \in \mathbf{R} $) on the vector space $\mathbf{R}^{N_{h} N_{k} }$ by
\begin{align}\label{Asnorm}
\left(\Vv_{k} , \Vw_{k}  \right)_{\bbA_{k} , s}  := \Vv_{k} ^{\top} \bbA_{k}^{s} \Vw_{k} , \quad
\|\Vv_{k} \|_{\bbA_{k} , s} ^{2}  := \left(\Vv_{k} , \Vv_{k} \right)_{\bbA_{k} ,s}.
\end{align}
By \eqref{Asnorm} and Lemma \ref{lem2.1}, we observe that
$\|\Vv_{k}\|_{\bbA_{k} ,1} ^{2} = \Ca_{k} (v_{k} , v_{k} ) \approx \N{v_{k}}_{1,k,h}^{2} . $ 

\subsection{Prolongation and restriction operators}
\label{sec3.1}
By \eqref{basis}, it is easy to see that
\begin{align}\label{L2-orth}
  ({\Psi^{\sigmabf}_T}, { \Psi^{\sigmabf'}_{T} })_{T}
= \prod_{i=1}^d \frac{\delta_{\sigma_i\sigma'_i} h}{2\sigma_i+1} \quad \forall \, \Psi_T^{\sigmabf} , \Psi_{T} ^{\sigmabf'}  \in U_{k} (\Ct_{h} ),
\end{align}
where $\delta_{\sigma_i\sigma'_i}$ is Kronecker symbol. 
Without confusions, we extend $\Psi^{\sigmabf}_T$ by zero to the exterior of $T$.
Combining \eqref{DGspace} with above identity, we find that  
\[
  U_{k} (\Ct_{h} ) = \Span{\Psi_{T} ^{\sigmabf} : \N{\sigmabf}_{\infty} \le k, \forall \, \sigmabf \in \mathbf{N}^{d}, \forall \, T \in \Ct_{h} }.
\]
For any level $k = 2, 3,\dots, p$, the prolongation operator connecting the space $U_{k-1} (\Ct_{h} )$ to $U_{k} (\Ct_{h} )$
is denoted by $\Cp_{k-1} ^{k} : U_{k-1} (\Ct_{h} ) \to U_{k} (\Ct_{h} )$. 
It is easy to see that
\[
  U_{k} (\Ct_{h} ) = U_{k-1}(\Ct_{h} ) \bigoplus \Span{\Psi_T^{\sigmabf} : \N{\sigmabf}_{\infty} = k, \forall \, \sigmabf \in \mathbf{N}^{d}, \forall \, T \in \Ct_{h} },
\] 
where $\bigoplus$ stands for the direct sum decomposition.
The above identity implies that the prolongation operator is just a natural injection. 
The restriction operator $\Cr_{k} ^{k-1}: U_{k} (\Ct_{h} ) \to U_{k-1} (\Ct_{h} )  $ is the adjoint of $\Cp_{k-1} ^{k} $ with respect to the $\Ltwo$ inner product, i.e.,
\[
  \SP{v_{k-1} }{\Cr_{k} ^{k-1} w_{k}}_{\Omega} := \SP{\Cp_{k-1} ^{k} v_{k-1} }{w_{k} }_{\Omega} \quad \forall \, v_{k-1} \in U_{k-1} (\Ct_{h} ), w_{k} \in U_{k} (\Ct_{h} ).  
\] 
The matrix representations of $\Cp_{k-1} ^{k}$ and $\Cr_{k} ^{k-1}$ 
are denoted by $\bbP_{k-1}^{k} , \bbR_{k} ^{k-1}$.
The $\Ltwo$-orthogonal property \eqref{L2-orth} implies that 
the entries of $\bbP_{k-1}^{k}$ are 0 or 1.
Moreover, $\bbR_{k} ^{k-1} = (\bbP_{k-1} ^{k} )^{\top}$.

\subsection{Smoothers}
\label{sec3.2}
For any level $k=1,2,\dots, p$, let $\bbD_{k} $ denote the diagonal matrix of $\bbA_{k} $ and $\bbI_{k} $ denote the identity matrix.
The classical algorithm of damped Jacobi smoother for $\bbA_{k} \Vz_{k}  = \Vf_{k}$ reads 
\begin{align}\label{smoothing_iteration}
\Vz_k^i = \Vz_k^{i-1} + \omega_{k} \bbD_{k} ^{-1} \left(\Vf_k - \bbA_k \Vz_k^{i-1} \right),
\end{align}
where $\Vz_k^0$ is the initial guess.
It is well-known that $\Vz_k^i$ converges for any initial guess if and only if $\omega_{k}  < 1 / \rho(\bbD_{k} ^{-1} \bbA_{k} )$. 
The matrix $ \left(\bbI_{k} - \omega_{k}  \bbD_{k}^{-1} \bbA_{k} \right)^{m} $ is called an iterative matrix for the damped Jacobi smoother \eqref{smoothing_iteration} with $m$ smoothing steps, i.e.,
\[
\Vz_{k}^{m} - \Vz_{k} = \left(\bbI_{k} - \omega_{k}  \bbD_{k} ^{-1} \bbA_{k} \right)^{m} \left(\Vz_{k} ^{0} - \Vz_{k} \right).
\]

Given an initial guess $\Vz_{k} ^{0} \in U_{k} (\Ct_{h} )$, we consider the fourth-kind Chebyshev polynomial iterative method (cf. \cite{Lottes2023})
\begin{align} \label{4th_Chebyshev}
\begin{cases}
\Vw_{k} ^{0} = \bm{0},\\
\Vw_k^{i} = \dfrac{2i-3}{2i+1}\Vw_k^{i-1} 
          + \dfrac{8i-4}{2i+1} \omega_{k} \bbD_{k}^{-1}
          \left(\Vf_k - \bbA_k \Vz_k^{i-1} \right), \\
\Vz_k^{i} = \Vz_k^{i-1} + \Vw_k^{i},\\ 
\end{cases}
\end{align}
  By $\Vf_{k} = \bbA_{k} \Vz_{k} $ and \eqref{4th_Chebyshev}, 
  for $i=1$ we get 
\[
  \Vz_{k}^{1} - \Vz_{k}  = 
  \big(\bbI_{k}  - \frac{ 4 }{ 3 } \omega_{k} \bbD_{k} ^{-1}  \bbA_{k} \big)
  \left(\Vz_{k} ^{0} - \Vz_{k} \right),
\]
  and  for $i\ge 2$ we have
\begin{align*}
\Vz_k^i - \Vz_k
&= \Vz_k^{i-1} - \Vz_k + \Vw_{k} ^{i} \\
&= \Vz_k^{i-1} - \Vz_k + \frac{2i-3}{2i+1}\Vw_k^{i-1} + \frac{8i-4}{2i+1} \omega_k \bbD_k^{-1} \bbA_k (\Vz_k - \Vz_k^{i-1})  \\
&= \big(\bbI_k - \frac{8i-4}{2i+1} \omega_k \bbD_k^{-1} \bbA_k \big) 
(\Vz_k^{i-1} - \Vz_k) 
+ \frac{2i-3}{2i+1}(\Vz^{i-1}_k - \Vz^{i-2}_k) \\
&= \frac{4i-2}{2i+1}\left(\bbI_k - 2 \omega_k \bbD_k^{-1} \bbA_k \right) (\Vz_k^{i-1} - \Vz_k)
- \frac{2i-3}{2i+1} (\Vz_k^{i-2} - \Vz_k).
\end{align*}
Then the iterative matrix of this polynomial smoother is a matrix polynomial 
$P_m(\omega_{k}  \bbD_{k}^{-1} \bbA_{k}  )$ of degree $m$, i.e. 
\[
\Vz_{k}^{m} - \Vz_{k} = P_m( \omega_{k}  \bbD_{k} ^{-1} \bbA_{k} ) \left(\Vz_{k} ^{0} - \Vz_{k} \right),
\] 
where $P_m(x)$ are defined recursively as follows 
\begin{equation} \label{polynomial_pm}
P_i(x) = \frac{4i-2}{2i+1}(1-2x)P_{i-1}(x) - \frac{2i-3}{2i+1}P_{i-2}(x) \quad i=1,2\cdots m,\\
\end{equation}
with $P_0(x) = 1$ and $P_1(x) = 1-\frac{4}{3}x$.
What's more, 
the supremums of polynomials $\sqrt{x} \left| P_m(x) \right|$ and $|P_{m} (x)|$ for $x\in [0,1]$ 
is given by 
\begin{equation}\label{poly_smoother1}
\sup\limits_{x\in [0,1] } \sqrt{x} \left|P_m(x) \right| = \frac{1}{2m+1},\quad 
\sup\limits_{x\in [0,1] } \left|P_m(x) \right| = 1 \quad \forall \, m\in \mathbf{N}^+.
\end{equation}
Later, we will use the polynomial smoother \eqref{4th_Chebyshev}
for $p$-MG method and show that it gives better smoothing property than the 
damped Jacobi smoother \eqref{smoothing_iteration}.

\begin{remark}[Complexity and storage of smoothers]
Comparing \eqref{smoothing_iteration} and \eqref{4th_Chebyshev}, 
we find that both smoothing algorithms require $m$ times of matrix-vector multiplications 
$\bbA_k \Vz_k^{i-1}, i=1,2,\cdots, m$. 
These are main computational quantities.
The polynomial smoother requires an additional storage for $\Vw_k^{i-1}$
from the previous iteration step.
This overhead is negligible by using dynamic programming. 
\end{remark}

\subsection{$W$-cycle algorithm}
\label{sec3.3}
\begin{algorithm}[!ht]
\caption{$W$-cycle algorithm}\label{Algorithm_k-MG}
\textbf{function} $W\text{cycle} (k, \bbA_k,\Vf_k,\Vz^0_{k}, m)$ 
\begin{algorithmic}[1]
    \IF{$k = 1$}
         \RETURN $\Vz = \bbA_{1} ^{-1} \Vf_{1}  $.
    \ELSE
        \FOR{$i = 1$ \TO $m$}
        \STATE $\Vz_k^{i} = \text{smoother}(\Vz_k^{i-1}, \bbA_k, \Vf_k)$.  
        \hfill \texttt{/*Pre-smoothing by \eqref{4th_Chebyshev}*/}
        \ENDFOR 
        \STATE 
        \STATE $\Vr_{k-1} = \bbR^{k-1}_{k}(\Vf_k - \bbA_k \Vz_k^{m})$. 
        \hfill \texttt{/*Restriction of the residual*/}
        \STATE $\Ve_{k-1} = W\text{cycle}(k-1, \bbA_{k-1}, \Vr_{k-1}, \bm{0}, m)$. 
        \hfill \texttt{/*Correction*/}
        \STATE $\Ve_{k-1} = W\text{cycle}(k-1, \bbA_{k-1}, \Vr_{k-1}, \Ve_{k-1}, m)$.
        \hfill \texttt{/*Correction*/}
        \STATE $\Vz_k^{m+1} = \Vz^m_k + \bbP_{k-1}^{k} \Ve_{k-1}$. 
        \hfill \texttt{/*Prolongation of the correction*/} 
        \STATE
        \FOR{$i = 1$ \TO $m$}
            \STATE $\Vz_k^{m+i+1} = \text{smoother}(\Vz_k^{m+i}, \bbA_k, \Vf_k)$. \hfill \texttt{/*Post-smoothing by 
            \eqref{4th_Chebyshev}*/}
        \ENDFOR
    \ENDIF
    \RETURN $\Vz \leftarrow \Vz_k^{2m+1}$
\end{algorithmic}
\end{algorithm}

The $W$-cycle algorithm of $p$-MG method is given by Algorithm \ref{Algorithm_k-MG}.
The error propagation matrix of Algorithm \ref{Algorithm_k-MG} is given by
\[
\bbE^{\rm W} _{k,m} \left(\Vz_{k}  - \Vz_{k} ^{0} \right) :=  
\Vz_{k}  - W\text{cycle} \left(k, \bbA_{k} , \Vf_{k}, \Vz_{k} ^{0} , m \right),
\] 
where $\Vz^{0}_k$ is the initial guess and $\Vz_{k} $ is the solution of $\bbA_{k} \Vz_{k}  = \Vf_{k} $. 
Let $\mathbb{G}_{k,m}^{\rm pre}, \mathbb{G}_{k,m}^{\rm post}$ 
denote the pre- and post-smoothing smoothing matrices, respectively. 
By Section \ref{sec3.2}, we have 
$\bbG_{k,m}^{\rm pre} = \bbG_{k,m}^{\rm post} = P_{m} \left(\omega_{k}  \bbD_{k} ^{-1} \bbA_{k} \right)$.
According to \cite{Elman2014}, we have the following recursive relation 
of the error propagation matrix 
\begin{align}\label{WError}
\bbE_{k,m}^{\rm W}
= \; & \bbG_{k,m}^{\rm post}
\big[\bbA^{-1}_k-\bbP_{k-1}^{k}
\big(\mathbb{I}_{k-1} - \left(\mathbb{E}^{\rm W}_{k-1,m}\right)^{2}  \big)
\bbA_{k-1}^{-1}\bbR_{k}^{k-1}\big]
\bbA_k \bbG_{k,m}^{\rm pre} \nonumber \\
= \; & 
\mathbb{E}^{\rm TL}_{k,m} + 
\bbG_{k,m}^{\rm post}\bbP_{k-1}^{k} \left(\mathbb{E}^{ \rm W}_{k-1,m} \right)^{2}
\bbA_{k-1}^{-1}\bbR_{k}^{k-1}\bbA_k \bbG_{k,m}^{\rm pre} ,
\end{align}
where $\mathbb{E}^{\rm W}_{k-1,m}$ is the error propagation matrix at $(k-1)$-level,
and $ \mathbb{E}^{\rm TL}_{k,m} $ is the error propagation matrix of the two-level algorithm, i.e.,
\begin{equation}\label{Twoerror}
\mathbb{E}^{\rm TL}_{k,m}
=\mathbb{G}_{k,m}^{\rm post}\left(\bbA^{-1}_k-\bbP_{k-1}^{k}\bbA_{k-1}^{-1}\bbR_{k}^{k-1}\right)\bbA_k \mathbb{G}_{k,m}^{\rm pre} .
\end{equation}

\section{Convergence analysis of $p$-MG algorithm}
\label{sec4}
In this section, we aim to prove the convergence of the $p$-MG algorithm. 
The proof consists of two steps: 
\begin{enumerate}
\item establishing the convergence of the two-level algorithm;
\item extending the theory to the $W$-cycle algorithm.
\end{enumerate} 
The convergence of the $p$-MG method relies primarily on two key aspects: the approximation property and the smoothing property.
We will present these studies first in order to provide a comprehensive convergence theory.

\subsection{Convergence of two-level algorithm}
\label{sec4.1}
The error propagation matrix of two-level algorithm 
$ \mathbb{E}^{\rm TL}_{k,m} $  
is the product of three error-damping matrices, 
coming from the pre-smoothing step, the correction step, 
and the post-smoothing step, respectively. 
The two-level algorithm converges if the contraction inequality
\[
  \N{\mathbb{E}^{\rm TL} _{k,m} \Vg_{k} }_{\bbA_k, 1} \leq C(h,k,m) \N{\Vg_{k} }_{\bbA_k, 1} 
  \quad \forall \, \Vg_{k} \in \mathbf{R}^{N_{h} N_{k} } ,
\]
holds with $ C(h,k,m) < 1 $. 
We will analyze how $ C(h,k,m) $ depends on the mesh size $ h $, 
the polynomial degree $ k $, and the number of smoothing steps
$ m $ in the two-level algorithm, 
and focus on the estimate of the error propagation operator. 

\subsubsection{Approximation property}
\label{sec4.1.1}
The mass matrix for the broken polynomial space $U_k(\Ct_h)$ is defined by 
\[
\Vv_{k} ^{\top} \bbM_{k} \Vw_{k} := \SP{v_{k}}{ w_{k} }_{\Omega} \quad \forall \, v_{k} ,w_{k} \in U_{k} (\Ct_{h}).
\] 
The $\Ltwo$-orthogonality \eqref{L2-orth} implies that 
the mass matrix $\bbM_{k} $ is diagonal, i.e. 
\[
\bbM_{k} = \diag\{ \NLtwo[T]{\Psi_{T} ^{\sigmabf} }^{2} : \N{\sigmabf}_{\infty} \le k, T \in \Ct_{h} \}.
\]
What's more, we have 
$\N{\bbM_{k} ^{\frac{1}{ 2 } } \Vv_{k} }_{\bbA_{k} ,0} ^{2} = \Vv_{k}^{\top}  \bbM_{k} \Vv_{k} = \NLtwo{v_{k} }^{2}$, for any $v_{k} \in U_{k} (\Ct_{h} ).$
Inspired by \cite[Lemma 2.8, Theorem 2.9]{Elman2014},
we obtain the approximation property. 
\begin{lemma}[Approximation property]\label{Lemma_approximation}
For any level $k=2,\cdots, p$,
there holds
\begin{equation}\label{approximation}
	\N{ \left(\bbA^{-1}_{k} - \bbP_{k-1}^{k}\bbA^{-1}_{k-1} \bbR_{k}^{k-1} \right) \Vg_k }_{\bbA_{k} ,1} 
	\lesssim \frac{h}{ k }  \N{ \mathbb{M}_k^{-\frac{1}{2}} \Vg_k }_{\bbA_{k} , 0} \quad
  \forall \, \Vg_{k} \in \mathbf{R}^{N_{h} N_{k} }
  \end{equation}
\end{lemma}

\begin{proof}
Suppose the components of $\Vg_{k} $ are given by the coefficients $\{ g_{T} ^{\sigmabf} : \N{\sigmabf}_{\infty} \le k, T\in \Ct_{h} \}$.
Define low-order vector $\Vg_{k-1} := \bbR_{k}^{k-1} \Vg_k \in \mathbf{R}^{N_{h} N_{k-1} }$ and piecewise polynomial 
\[
\tilde{g}_{k} (\Bx) :=  \sum_{T\in\Ct_h}\sum_{\N{\sigmabf}\leq k} \frac{ g_{T} ^{\sigmabf} }{ \NLtwo[T]{\Psi_{T} ^{\sigmabf} }^{2}   }   \Psi_T^{\sigmabf} (\Bx).
\]
It is easy to verify 
\begin{align}\label{eq-1}
\NLtwo[\Omega]{\tilde{g}_{k}}^2 =  \Vg_{k} ^{\top} \bbM_k^{-1} {\Vg_{k} },\quad
\SP{\tilde{g}_{k}}{\Psi_{T} ^{\sigmabf}}_{T} = g_{T} ^{\sigmabf} \quad \forall \, \Psi_{T} ^{\sigmabf} \in U_{k} (\Ct_{h} ).
\end{align}
Define $\Vu_{k}  = \bbA_{k} ^{-1} \Vg_{k} $ and 
$\Vu_{k-1}  = \bbA_{k-1} ^{-1} \bbR_{k} ^{k-1} \Vg_{k} $
and let $u_{k} \in U_{k} (\Ct_{h} ), u_{k-1} \in U_{k-1} (\Ct_{h} )$ be the functions 
corresponding to $\Vu_k, \Vu_{k-1}$, respectively, according to \eqref{LinearSystem}.
By \eqref{Asnorm} and the natural injection $\Cp_{k-1} ^{k} $, we find that
\begin{align}\label{approx-1}
\begin{split}
\N{ \left(\bbA^{-1}_{k} - \bbP_{k-1}^{k}\bbA^{-1}_{k-1} \bbR_{k}^{k-1} \right) \Vg_k }_{\bbA_k,1}
& = \N{ \Vu_{k} - \bbP_{k-1}^{k}\Vu_{k-1} }_{\bbA_k,1} \\
& \lesssim \N{u_k - \mathcal{P}_{k-1}^{k}u_{k-1} }_{1,k,h} \\
& = \N{u_{k} - u_{k-1}}_{1,k,h}.
\end{split}
\end{align}

It is easy to see that $u_k$ and $u_{k-1}$ are the discrete approximates of the solution to the boundary-value problem
\begin{equation}
-\Delta u = \tilde{g}_k \quad \text{in}\;\; \Omega, \quad
u=0\quad \text{on}\;\;  \partial \Omega,
\end{equation}
in $U_k(\Ct_h)$ and $U_{k-1}(\Ct_h)$, respectively.
By triangle inequality and the optimal error estimate \eqref{error_estimate_1}, we have
\begin{align}\label{ineq-1}
\N{u_{k} - u_{k-1} }_{1,k,h} 
\lesssim \N{u_k - u}_{1,k,h} + \N{u - u_{k-1}}_{1,k,h}
\lesssim \frac{ h }{ k }  \N{u}_{H^2(\Ct_{h})},
\end{align}
Combine \eqref{eq-1}, \eqref{ineq-1} and the regularity property $\N{u}_{H^2(\Ct_{h} )} \lesssim \NLtwo{\tilde{g}_{k} }$ yields 
\begin{align}\label{ineq-2}
  \N{u_{k} - u_{k-1} }_{1,k,h} \lesssim 
  \frac{h}{ k } \N{ \bbM_{k} ^{-\frac{1}{ 2 } } \Vg_{k} }_{\bbA_{k} ,0} .
\end{align}
The proof is complete by \eqref{approx-1} and \eqref{ineq-2}.
\end{proof}

\subsubsection{Smoothing properties}
\label{sec4.1.2}
We consider the damping matrix $\omega_k \mathbb{D}^{-1}_k$,
where $\omega_k$ is the damping parameter satisfying $0 < \omega_k \rho(\mathbb{D}^{-1}_k \bbA_k) \leq 1$.
For the estimation of spectral radius $\rho(\bbD_{k} ^{-1} \bbA_{k})$,
is given by the following lemma. 

\begin{lemma}\label{Lemma_omega-k}
For any level $k=2,\dots,p$, there holds $\rho(\bbD_{k} ^{-1} \bbA_{k} ) \lesssim k+1$ . 
\end{lemma}

The proof of this Lemma \ref{Lemma_omega-k} is inspired by \cite[Proposition1.10]{Vassilev2008}, 
and the details are given in Appendix \ref{app}. 
We will take {$\omega_k \lesssim {1}/{(k+1)}$} 
in the remaining analysis.
Before studying the smoothing properties of the polynomial smoother,
we present two useful lemmas to aid the analysis. 

\begin{lemma}\label{Lemma_exchange}
For any level $k=2,\dots,p$, define $\bbC_{k} := \bbD_{k} ^{-\frac{1}{2}}\bbA_{k} \bbD_{k} ^{-\frac{1}{2}}$. 
Then 
\begin{align}\label{exchange}
\bbA_{k}  P_m(\omega_{k} \bbD_{k} ^{-1}\bbA_{k} ) 
= \bbD_{k} ^{\frac{1}{2}}\bbC_{k}  P_m(\omega_{k} \bbC_{k}  ) \bbD_{k} ^{\frac{1}{2}}.
\end{align}
\end{lemma}

\begin{proof}
Since $\bbA_{k} $ is a symmetric and positive definite (SPD) matrix, 
$\bbD_{k}$ and $\bbD_{k}^{\pm\frac{1}{2}}$ exist and are both SPD.
First, we show 
\begin{align}\label{matrix_DC}
\bbA_{k} (\bbD_{k} ^{-1}\bbA_{k})^m 
= \bbD_{k} ^{\frac{1}{2}} \bbC_{k} ^{m+1} \bbD_{k} ^{\frac{1}{2}}.
\end{align}
For $m=1$, clearly we have 
$\bbA_{k} \bbD_{k} ^{-1} \bbA_{k}  
= \bbD_{k} ^{\frac{1}{2}} \bbC_{k} ^2 \bbD_{k} ^{\frac{1}{2}}$. 
Suppose \eqref{matrix_DC} holds for $m-1$, then we have 
$\bbA_{k} \left( \bbD_{k} ^{-1}\bbA _{k} \right)^m 
= \bbD_{k} ^{\frac{1}{ 2 } } \bbC_{k} ^{m} \bbD_{k} ^{\frac{1}{ 2 } } (\bbD_{k} ^{-1} \bbA_{k} )
= \bbD_{k} ^{\frac{1}{2}} \bbC_{k} ^{m+1} \bbD_{k} ^{\frac{1}{2}}$, 
which yields \eqref{matrix_DC} for $m$.
The proof is finished by induction. 
\end{proof}

\begin{lemma}\label{Lemma_invDM}
For any level $k=2,\dots,p$,
there holds
\begin{align}\label{theta}
\rho(\bbM^{-1}_k\bbD_k) \lesssim \frac{k^3}{ h^2}.
\end{align}
\end{lemma}

\begin{proof}
By \eqref{basis} straightforward calculations, we find that 
\begin{align}\label{diag-e}
\Ca_{k} (\Psi_{T} ^{\sigmabf}, \Psi_{T} ^{\sigmabf} ) = 
\frac{2\alpha_0 k^2}{h^2} 
\big(\prod_{i=1} ^{d} \frac{ h }{ 2\sigma_i  + 1} \big) 
\sum_{i=1} ^{d} \left(2 \sigma_i + 1\right) \quad 
\forall \, \Psi_{T} ^{\sigmabf} \in U_{k} (\Ct_{h} ).
\end{align} 
The details are given in Appendix \ref{app}.
It is clear that $\bbM_{k} ^{-1} \bbD_{k} $ is diagonal, i.e. 
\[
 \bbM_{k} ^{-1} \bbD_{k} = 
 \diag \left\{  
 \SP{\Psi_{T} ^{\sigmabf} }{\Psi_{T} ^{\sigmabf} }_{T} ^{-1} 
 \Ca_{k} \left(\Psi_{T} ^{\sigmabf} , \Psi_{T} ^{\sigmabf} \right):
 \N{\sigmabf}_\infty \le k, T \in \Ct_{h} 
 \right\}.
\] 
Combine the $\Ltwo$-orthogonality \eqref{L2-orth} and \eqref{diag-e}. 
The diagonal entries of $\bbM_k^{-1}\bbD_k$ are given by
\begin{align}\label{invM_D}
\SP{\Psi_{T} ^{\sigmabf} }{\Psi_{T} ^{\sigmabf} }_{T} ^{-1} 
\Ca_{k} \left(\Psi_{T} ^{\sigmabf} ,\Psi_{T} ^{\sigmabf} \right)
=   \frac{ 2 \alpha_0 k^{2}   }{ h^{2}  }  \sum_{i=1} ^{d} \left(2 \sigma_i + 1\right)   
\quad \forall \, \Psi_{T} ^{\sigmabf} \in U_{k} (\Ct_{h} ).
\end{align}
Since $\sigma_i \leq k$ for $i=1, \dots, d$, we get \eqref{theta}. 
\end{proof}

Based on these lemmas, we have the following two smoothing properties for the polynomial smoother $\bbG_{k,m} ^{\rm pre} = \bbG_{k,m} ^{\rm post} = P_m(\omega_k \mathbb{D}^{-1}_k \bbA_k )$.
\begin{lemma}[pre-smoothing property]\label{Lemma_pre-smoothing}
Suppose $\omega_k \lesssim (k+1)^{-1} $. Then 
\begin{eqnarray}\label{pre-smoothing}
\N{ \bbM_{k} ^{-\frac{1}{2}} \bbA_k P_{m}\left(\omega_k \bbD_k^{-1}\bbA_k \right) \Vg_{k} }_{\bbA_k,0} \lesssim \frac{ k^{2}  }{ h m  }  \N{\Vg_{k} }_{\bbA_k,1} \quad 
\forall \, \Vg_{k}  \in \mathbf{R}^{N_{h} N_{k} } .
\end{eqnarray}
\end{lemma}

\begin{proof}
Since $\bbM_{k} $ is diagonal, 
we have $\bbD_{k}^{\frac{1}{ 2 } } \bbM_{k} ^{-1} \bbD_{k} ^{\frac{1}{ 2 } } = \bbM_{k}^{-1}  \bbD_{k} $.
By Lemma \ref{Lemma_exchange}, we have
\begin{align}\label{ineq-3}
\N{ \bbM^{-\frac{1}{2}}_k \bbA_k P_{m}\left(\omega_k \bbD_k^{-1} \bbA_k \right) \Vg_{k} }_{\bbA_k,0}^2
\leq \rho\left(\bbM^{-1}_k\bbD_k \right) 
\N{ \bbC_k P_{m}\left(\omega_k \bbC_k \right) \bbD^{\frac{1}{2}}_k\Vg_{k} }^2_{\bbA_k,0} .
\end{align}

Suppose $\{\lambda_i\}_{i=1}^{N_hN_k}, \{\bm{\varphi}_i\}_{i=1}^{N_hN_k}$  are eigenvalues and orthogonal eigenvectors of $\bbC_k$, satisfying 
\begin{align}\label{eignvs}
  \bbC_k \bm{\varphi}_i = \lambda_i \bm{\varphi}_i, \quad (\bm{\varphi}_i, \bm{\varphi}_j)_{\bbA_{k} ,0}  = \delta_{ij} \quad \forall \, i , j \in \{1,2,\dots,N_{h} N_{k}\},
\end{align}
On the one hand, 
$\bbD_{k} ^{-\frac{1}{ 2 } } \bbC_{k} \bbD_{k} ^{\frac{1}{ 2 } } = \bbD_{k} ^{-1} \bbA_{k} $, 
implies $\omega_{k} \rho (\bbC_{k} ) = \omega_{k} \rho(\bbD_{k} ^{-1} \bbA_{k} ) \leq 1 $.
On the other hand, 
$\Vw_{k} ^{\top} \bbC_{k} \Vw_{k} = (\bbD_{k} ^{ - \frac{1}{ 2 }} \Vw_{k} )^{\top}  \bbA_{k} (\bbD_{k} ^{- \frac{1}{ 2 } } \Vw_{k} ) >0$,
for any $\Vw_{k} \neq \mathbf{0}.$
Therefore, we have $\omega_{k} \lambda_{i} \in (0,1], \; 1\leq i\leq N_h N_k$. 
Write $\Vv_{k}  = \bbD^{\frac{1}{2}}_k\Vg_{k}  $, 
$\Vv_{k} = \sum_{i=1}^{N_hN_k}v_i \bm{\varphi}_i$,
according to the property of the polynomial $P_m(\omega_k \bbC_k) \bm{\varphi}_i = P_m(\omega_k \lambda_i) \bm{\varphi}_i$, 
we have 
$$
\bbC_k P_m(\omega_k \bbC_k) \Vv_k 
= \sum_{i=1} ^{N_{h} N_{k}} P_{m} (\omega_{k} \lambda_{i} ) v_{i} \bbC_k \varphibf_{i}
= \sum_{i=1} ^{N_{h} N_{k}} \lambda_{i} P_{m} (\omega_{k} \lambda_{i} ) v_{i} \varphibf_{i}. 
$$
Then, using the supremum \eqref{poly_smoother1} of the polynomial $\sqrt{x} P_{m} (x)$, we have 
\begin{align}\label{ineq-4}
\begin{split}
\N{ \bbC_k P_{m}\left(\omega_k \bbC_k \right) \Vv_{k} }^2_{\bbA_k,0} 
&=\bigg(\sum_{i=1} ^{N_{h} N_{k}} \lambda_{i} P_{m} (\omega_{k} \lambda_{i} ) v_{i} \varphibf_{i}  \bigg) ^{\top} 
\bigg(\bbC_{k} \sum_{i=1}^{N_{h} N_{k} }P_{m} (\omega_{k} \lambda_{i} ) v_{i}  \varphibf_{i}  \bigg)\\
&=\bigg(\sum_{i=1} ^{N_{h} N_{k}} \lambda_{i} P_{m} (\omega_{k} \lambda_{i} ) v_{i} \varphibf_{i}  \bigg) ^{\top} 
\bigg(\sum_{i=1}^{N_{h} N_{k} }P_{m} (\omega_{k} \lambda_{i} ) v_{i} \bbC_{k}  \varphibf_{i}  \bigg)\\
&\lesssim 
\bigg( \sup\limits_{\omega_{k} \lambda_i\in(0,1]} {\omega_{k} \lambda_i P^2_m(\omega_{k} \lambda_i)} \bigg)
\frac{\Vv_{k} ^{\top} \bbC_k \Vv_{k}}{ \omega_{k}  } \\
&= \frac{\Vg_{k} ^{\top} \bbA_{k} \Vg_{k} }{\omega_k (2m+1)^2} 
\lesssim \frac{k\N{\Vg_{k} }^2_{\bbA_k,1}}{m^2},
\end{split}
\end{align}
where we have used $\omega_k \lesssim ({k+1})^{-1} $ in the last inequality. 
The proof is finished by combing \eqref{ineq-3}, \eqref{ineq-4} and Lemma \ref{Lemma_invDM}. 
\end{proof}

\begin{lemma}[post-smoothing property]\label{Lemma_post-smoothing}
For any level $k=2,\dots, p$, we have
\begin{align}\label{post-smoothing}
\N{ P_m(\omega_k \mathbb{D}^{-1}_k\bbA_k) \Vg_{k} }_{\bbA_k, 1} 
\leq \N{ \Vg_{k}  }_{\bbA_k, 1}\quad
\forall \, \Vg_{k} \in \mathbf{R}^{N_{h} N_{k}}. 
\end{align}
\end{lemma}

\begin{proof}
Similar to the proof of Lemma \ref{Lemma_pre-smoothing}, 
assuming $\{\lambda_i\}_{i=1}^{N_hN_k}, \{\bm{\varphi}_i\}_{i=1}^{N_hN_k}$ are corresponding eigenvalues and orthogonal eigenvectors of $\bbD^{-1}_k\bbA_k$,
we have $P_m\left( \omega_k \bbD^{-1}_k\bbA_k \right) \bm{\varphi}_i 
= P_m(\omega_{k} \lambda_i) {\varphibf}_i$.
Write $\Vg_{k} $ as $\Vg_{k}  = \sum_{i=1}^{N_{h}N_{k}} g_i \bm{\varphi}_i$.
Then, we have 
\begin{align*}
\N{P_m(\omega_k \mathbb{D}^{-1}_k\bbA_k) \Vg_{k}}_{\bbA_k,1} 
&=\N{ \sum_{i=1}^{N_hN_k} g_i P_m\left(\omega_{k} \lambda_i \right) \bm{\varphi}_i }_{\bbA_k,1} \\
&\leq \bigg(\sup_{\omega_{k} \lambda_{i} \in (0,1]} P_{m} (\omega_{k} \lambda_{i} )\bigg)
\N{\sum_{i=1}^{N_hN_k} g_i \bm{\varphi}_i}_{\bbA_k,1} .
\end{align*}
By the supremum \eqref{poly_smoother1} of the polynomial $P_{m} (x)$, 
we get \eqref{post-smoothing}.
\end{proof}

\subsubsection{Convergence of the two-level algorithm}
\label{sec4.1.3}
Combining the approximation property and the smoothing properties, 
we have the following convergence rate for the two-level algorithm.
\begin{theorem}\label{Thm_two-level}
Suppose the number of pre-smoothing steps $m\gtrsim k$.
The two-level algorithm converges 
with rate independent of $h$ and $k$. 
\end{theorem}

\begin{proof}
By the approximation property \eqref{approximation} and the smoothing properties \eqref{pre-smoothing}, \eqref{post-smoothing}, 
the error propagation matrix $\mathbb{E}^{\rm TL}_{k,m}$ given by \eqref{Twoerror} satisfies, 
\begin{align}\label{Twoconv}
\begin{split}
\N{ \mathbb{E}^{\rm TL}_{k,m} \Vg_{k} }_{\bbA_k, 1}
&=\N{ \bbG_{k,m}^{\rm post}\left(\bbA^{-1}_k-\bbP_{k-1}^{k}\bbA_{k-1}^{-1}\bbR_{k}^{k-1}\right)\bbA_k \bbG_{k,m}^{\rm pre}\Vg_{k}  }_{\bbA_k, 1} \\
&\leq\N{ \left(\bbA^{-1}_k-\bbP_{k-1}^{k}\bbA_{k-1}^{-1}\bbR_{k}^{k-1}\right)\bbA_k P_m(\omega_k \mathbb{D}^{-1}_k\bbA_k) \Vg_{k}  }_{\bbA_k, 1} \\
&\leq C \frac{h}{k} \N{ \mathbb{M}_k^{-\frac{1}{2}}\bbA_k 
P_m(\omega_k \mathbb{D}^{-1}_k\bbA_k) 
\Vg_{k}  }_{\bbA_k, 0} \\
&\leq C \frac{h}{k} \frac{k^2}{hm} \N{ \Vg_{k}  }_{\bbA_k, 1} 
= C \frac{k}{m} \N{ \Vg_{k}  }_{\bbA_k, 1}, 
\end{split}
\end{align}
where $C > 0$ is a constant independent of $h$ and $k$.
Taking $m > Ck$ shows that 
$\N{\bbE_{k,m} ^{\rm TL} \Vg_{k}}_{\bbA_{k},1} / \N{\Vg_{k} }_{\bbA_{k} ,1} < 1 $.
Consequently, the two-level algorithm converges with rate independent of both $h$ and $k$. 
\end{proof}

\subsection{Convergence of the $W$-cycle algorithm}
\label{sec4.2}
In this subsection, we study the $W$-cycle algorithm. 
First, we need some estimates about the stability of the prolongation and restriction operators.
\begin{lemma}\label{Lemma_stability}
For any level $k=2,\dots,p$, 
there holds
\begin{align}
\N{\bbP_{k-1}^{k}\Vv_{k-1}  }_{\bbA_{k},1} 
&\lesssim \N{\Vv_{k-1} }_{\bbA_{k-1},1}
\quad \forall \, \Vv_{k-1} \in \mathbf{R}^{N_{h} N_{k-1} }, \label{stab_1}\\ 
\N{\bbA_{k-1}^{-1}\bbR_{k}^{k-1}\bbA_k\Vw_{k}  }_{\bbA_{k-1},1}
&\lesssim \N{\Vw_{k} }_{\bbA_{k},1}
\quad \quad \;\;\,\forall \, \Vw_{k} \in \mathbf{R}^{N_{h} N_{k} } . 
\label{stab_2}
\end{align}
\end{lemma}

\begin{proof}
Let $\Vv_{k-1} $ be the vector associated with $v_{k-1} \in U_{k-1} (\Ct_{h} )$,
according to \eqref{basisexp} and satisfying 
$\N{\Vv_{k-1} }_{\bbA_{k-1} ,1}^{2} = \Ca_{k-1} (v_{k-1},v_{k-1}) $.
By 
$\Cp_{k-1}^{k} v_{k-1}  = v_{k-1} $ and Lemma \ref{lem2.1}, 
we have
\[
  \N{\bbP_{k-1} ^{k} \Vv_{k-1} }_{\bbA_{k} ,1} ^{2} 
  = \Ca_{k} (\Cp_{k-1} ^{k} v_{k-1} , \Cp_{k-1} ^{k} v_{k-1} )
  \lesssim \N{\Cp_{k-1} ^{k} v_{k-1} }_{1,k,h}^{2} 
  \lesssim \N{v_{k-1} }_{1,k-1,h} ^{2} .
\] 
So \eqref{stab_1} easily follows from the coercivity of $\Ca_k$ in Lemma \ref{lem2.1}.

By Schwarz's inequality, we have 
\begin{align*}
\N{\bbA_{k-1}^{-1}\bbR_{k}^{k-1}\bbA_k\Vw_{k}}_{\bbA_{k-1},1}
= & \max\limits_{\Vu_{k-1} \in \mathbf{R}^{N_{h} N_{k-1}} \slash \{\mathbf{0}\}}
\frac{ \left( \bbA_{k-1}^{-1}\bbR_{k}^{k-1}\bbA_k\Vw_{k} , \Vu_{k-1} \right)_{\bbA_{k-1} ,1}  }
{\N{\Vu_{k-1} }_{\bbA_{k-1}, 1} } \\
= & \max\limits_{\Vu_{k-1} \in \mathbf{R}^{N_{h} N_{k-1}} \slash \{\mathbf{0}\}}
\frac{\Vw_{k} ^{\top} \bbA_{k} \bbP_{k-1} ^{k} \Vu_{k-1}  } 
{\N{\Vu_{k-1} }_{\bbA_{k-1}, 1} } \\
\leq & \max\limits_{\Vu_{k-1} \in \mathbf{R}^{N_{h} N_{k-1}} \slash \{\mathbf{0}\}} 
\frac{\N{\bbP_{k-1}^{k}\Vu_{k-1} }_{\bbA_k, 1} \N{\Vw_{k} }_{\bbA_k, 1} }
{\N{\Vu_{k-1} }_{\bbA_{k-1}, 1} }. 
\end{align*}
Combining \eqref{stab_1} and the above inequality, we get \eqref{stab_2}.
\end{proof}

Now we are ready to prove the main result of this paper concerning the convergence rate of the $W$-cycle algorithm of $p$-MG method. 

\begin{theorem}\label{Thm_W-cycle}
There exists a constant $C_0 > 0$ independent of $h$ and $k$, 
such that, when $m > C_0 k$, the $W$-cycle algorithm of $k$-level $p$-MG converges 
with rate $\delta\in (0, 1)$ independent of $h$ and $k$, i.e. 
\begin{equation}\label{W-cycle}
\N{\bbE_{k,m}^{\rm W}\Vg_{k} }_{\bbA_{k},1} 
< \delta \N{\Vg_{k} }_{\bbA_{k},1} \quad
\forall \, \Vg_{k} \in \mathbf{R}^{N_{h} N_{k} } .
\end{equation}
\end{theorem}

\begin{proof}
For $k=2$, the W-cycle algorithm turns into the two-level algorithm, 
so \eqref{W-cycle} holds naturally.
For $k > 2$, suppose \eqref{W-cycle} holds for the ($k-1$)-level algorithm. 
Next we use induction to prove the conclusion \eqref{W-cycle} for the $k$-level algorithm. 
Use \eqref{WError}, we have 
\begin{align}\label{ineq-5}
\N{\bbE_{k,m}^{\rm W}\Vg_{k}  }_{\bbA_{k},1} 
\leq \N{\mathbb{E}^{\rm TL}_{k-1,m}\Vg_{k} }_{\bbA_{k},1}  + 
\N{\bbG_{k,m}^{\rm post}\bbP_{k-1}^{k} 
\left(\mathbb{E}^{\rm W}_{k-1,m}\right)^{2}
\bbA_{k-1}^{-1}\bbR_{k}^{k-1}\bbA_k \bbG_{k,m}^{\rm pre}\Vg_{k}  }_{\bbA_{k},1}. 
\end{align}
By inductive hypothesis, suppose we have
$$
\N{\bbE_{k-1,m}^{\rm W}\Vg_{k-1} }_{\bbA_{k-1},1} 
< \frac{k-1}{m} \N{\Vg_{k-1} }_{\bbA_{k-1},1} \quad
\forall \, \Vg_{k-1} \in \mathbf{R}^{N_{h} N_{k-1} } .
$$
Combine with the smoothing property \eqref{post-smoothing}, 
and the stabilities \eqref{stab_1}-\eqref{stab_2}, we have
\begin{equation}
\begin{aligned}\label{ineq-6}
&\N{\bbG_{k,m}^{\rm post}\bbP_{k-1}^{k} 
\left(\mathbb{E}^{\rm W}_{k-1,m}\right)^{2}
\bbA_{k-1}^{-1}\bbR_{k}^{k-1}\bbA_k \bbG_{k,m}^{\rm pre}\Vg_{k}  }_{\bbA_{k},1} \\
\lesssim &\N{ \left(\mathbb{E}^{\rm W}_{k-1,m}\right)^{2}
\bbA_{k-1}^{-1}\bbR_{k}^{k-1}\bbA_k \bbG_{k,m}^{\rm pre}\Vg_{k}  }_{\bbA_{k},1} \\
\lesssim &\frac{(k-1)^2}{m^2}\N{ \bbA_{k-1}^{-1}\bbR_{k}^{k-1}\bbA_k \bbG_{k,m}^{\rm pre}\Vg_{k}  }_{\bbA_{k},1} \\
\lesssim &\frac{ (k-1)^{2}  }{ m^{2}  }  \N{\Vg_{k}  }_{\bbA_{k},1} .
\end{aligned}
\end{equation}
Taking $m \gtrsim k$ and combining \eqref{Twoconv}, \eqref{ineq-6}, and \eqref{ineq-5}, we have
\begin{align*}
\N{\bbE_{k,m}^{\rm W}\Vg_{k}  }_{\bbA_{k},1} 
\lesssim \left( \frac{ k }{ m } + \frac{ (k-1)^{2}  }{ m^{2}  }  \right)\N{\Vg_{k} }_{\bbA_{k},1} 
\lesssim \frac{ k }{ m } \N{\Vg_{k} }_{\bbA_{k},1}.
\end{align*}
This completes the proof.
\end{proof}

\section{$p$-MG with the inherited bilinear form} 
\label{sec5}
Now, we consider the following matrices at different levels generated by the inherited bilinear form. 
For $k=1,\dots,p-1$, 
\[
  \Vv_{k}^{\top}  \tilde{\bbA}_{k} \Vw_{k} :=  \tilde{\Ca}_k(v_k,w_k), \quad 
  \tilde{\Ca}_{k} \left(v_{k} ,w_{k} \right) := \tilde\Ca_{k+1}({\Cp_{k}^{k+1}v_k, \Cp_{k}^{k+1}w_k})
  \quad \forall \, v_{k} ,w_{k} \in U_{k} (\Ct_{h}),
\]
where $\Vv_{k} ,\Vw_{k} $ are the vectors of unknowns of $v_{k} ,w_{k} $,  respectively, 
and $ \tilde\Ca_{p} (v_{p} ,w_{p} ) = \Ca_{p} (v_{p} ,w_{p} ) $ for any $v_{p} ,w_{p} \in U_{p} (\Ct_{h} )$.
Obviously, $\{\tilde{\bbA}_k\}_{k=1} ^{p-1} $ are all sub-matrices of ${\bbA}_p$.
This makes the algorithm very efficient. 
The $W$-cycle algorithm of $p$-MG with inherited bilinear forms is almost as same as Algorithm \ref{Algorithm_k-MG}, 
by replacing $\{ \bbA_{k} \}_{k=1} ^{p-1} $ with $\{\tilde{\bbA}_{k} \}_{k=1} ^{p-1} $.
In the meantime, we have the following modified continuity, coercivity, and error estimates.

\begin{lemma}
Suppose $\alpha_0$ is large enough. 
For any level $k = 1,\dots, p$, the inherited bilinear forms admit, 
\begin{align*}
\tilde{\Ca}_{k} (v_k,v_k) \gtrsim \N{v_k}^2_{1,p,h},\quad
\SN{\tilde{\Ca}_{k} (v_k,w_k)} \lesssim \N{v_k}_{1,p,h}\N{w_k}_{1,p,h}
\quad \forall \, v_{k} ,w_{k}  \in U_{k} (\Ct_{h} ).
\end{align*}
Whereas the error estimate bound modifies as follows:
\begin{align*}
  \N{u-\tilde u_{k} }_{1,p,h} \lesssim \frac{ h p}{ k^{2}  }  \N{u}_{H^{2} (\Ct_{h} )} , 
\end{align*}
where $\tilde u_{k} $ is the solution to the discrete problem,
\[
  \tilde{\Ca}_{k} (\tilde u_{k} ,v_{k} ) = (f, v_{k} ) \quad \forall \, v_{k} \in U_{k} (\Ct_{h} ).
\] 
\end{lemma}

Based on the above inherited bilinear forms, 
we can also prove the convergence of the two-level and $W$-cycle algorithms.
The details are omitted here. 

\begin{lemma}[approximation property]\label{Lemma_approximation_inherited}
For any level $k=2,\dots,p$, we have
\begin{align*}
\N{ \left( \tilde{\bbA}^{-1}_{k} - \bbP_{k-1}^{k} \tilde{\bbA}^{-1}_{k-1} \bbR_{k}^{k-1} \right) \Vg_k }_{\tilde{\bbA}_{k}  ,1} 
	\lesssim \frac{hp}{k^{2} }\N{ {\bbM}_k^{-\frac{1}{2}} \Vg_k }_{ \tilde{\bbA}_{k} , 0}
   \quad
  \forall \, \Vg_{k} \in \mathbf{R}^{N_{h} N_{k} }. 
\end{align*}
\end{lemma}

\begin{lemma}[smoothing properties]\label{Lemma_pre_smoothing_inherited}
Suppose $\tilde{\omega}_{k} \lesssim (p+1)^{-1} $, for any level $k=2,\dots, p$, we have
\begin{align*}
\N{ {\bbM}^{-\frac{1}{2}}_k \tilde{\bbA}_k 
P_{m}\left(\tilde{\omega}_k \tilde{\bbD}_k^{-1} \tilde{\bbA}_k \right) \Vg_{k} }_{ \tilde{\bbA}_k,0} 
& \lesssim  \frac{ kp }{ hm } \N{\Vg_{k} }_{\tilde{\bbA}_k,1}  \quad 
  \forall \, \Vg_{k} \in \mathbf{R}^{N_{h} N_{k} },\\
\N{ P_m( \tilde{\omega}_k \tilde{\bbD}^{-1}_k \tilde{\bbA}_k) \Vg_{k}  }_{ \tilde{\bbA}_k, 1} 
& \leq \; \N{ \Vg_{k}  }_{ \tilde{\bbA}_k, 1} \quad
\quad \;\, \forall \, \Vg_{k} \in \mathbf{R}^{N_{h} N_{k} },
\end{align*}
where $\tilde{\bbD}_{k} $ is the diagonal matrix of $\tilde{\bbA}_{k}$.
\end{lemma}

\begin{theorem}\label{Thm_two-level-inherited}
With the inherited bilinear forms, 
the two-level algorithm converges 
when the number of pre-smoothing steps $m \gtrsim p^{2} k^{-1} $
and the convergence rate is independent of $h$ and $k$. 
\end{theorem}

\begin{lemma}\label{Lemma_stability_inherited}
For any level $k=2,\dots, p$, we have
\begin{align*}
\N{\bbP_{k-1}^{k}\Vv_{k-1}  }_{\tilde{\bbA}_{k},1} 
&  = \N{\Vv_{k-1} }_{\tilde{\bbA}_{k-1},1}
\quad \forall \, \Vv_{k}  \in \mathbf{R}^{N_{h} N_{k}} 
,  \\
\N{\tilde{\bbA}_{k-1}^{-1}\bbR_{k}^{k-1} \tilde{\bbA}_k\Vw_{k}  }_{ \tilde{\bbA}_{k-1},1}
& \lesssim \N{\Vw_{k} }_{\tilde{\bbA}_{k},1}\quad  
\;\;\quad  \forall \, \Vw_{k} \in \mathbf{R}^{N_{h} N_{k} } .
\end{align*}
\end{lemma}

\begin{theorem}\label{Thm_W-cycle_inherited}
With the inherited bilinear forms, 
the W-cycle algorithm for $k$-level of $p$-MG method converges if 
$m \gtrsim {p^{2}} k^{-1} $, i.e.
\begin{align*}\label{W-cycle_inherited}
\N{ \tilde{\bbE}_{k,m}^{\rm W}\Vg_{k} }_{ \tilde{\bbA}_{k},1} 
\lesssim \frac{ p^{2}  }{k m }  \N{\Vg_{k}  }_{ \tilde{\bbA}_{k},1} 
\quad \forall \, \Vg_{k} \in \mathbf{R}^{N_{h} N_{k} }. 
\end{align*}
\end{theorem}

We will verify our $W$-cycle algorithms in the next section. 
It can be observed that
the convergence rates of $W$-cycle algorithms with the inherited and non-inherited forms 
are essentially the same.
However, in practical computations, 
the $W$-cycle algorithm with the inherited forms offers higher computational efficiency 
and requires less storage space
than the one using non-inherited bilinear forms,
as it does not need to build stiffness matrices on coarse levels.


\section{Numerical experiments}
\label{sec6}
In this section, we present some numerical tests 
to verify theoretical results and convergence rates of our $p$-MG algorithm.
One- and two-dimensional numerical experiments are performed on a workstation with a 2.9 GHz CPU (intel Xeon Gold 6326, 16 processors), and the codes are written by MATLAB language without parallel implementation.
Three-dimensional numerical tests are carried out on the LSSC-IV cluster at the State Key Laboratory of Scientific and Engineering Computing (LSEC),
Chinese Academy of Sciences, and the codes are implemented by adaptive finite element package ``Parallel Hierarchical Grid'' (PHG) \cite{Zhang2009}. 

In the numerical experiments, 
the $W$-cycle algorithm starts on the $p$-level
with non-inherited bilinear forms $\{\Ca_k\}_{k=1} ^{p} $ 
or inherited bilinear form $\{\tilde{\Ca }_{k}\}_{k=1} ^{p} $.
Two-level algorithm for $p$-level linear system $\bbA_{p} \Vu_{p} = \Vf_{p} $ 
uses directly solver for
the $(p-1)$-level linear system $\bbA_{p-1}\bm{e}_{p-1} = \Vr_{p-1}$ or $\tilde{\bbA}_{p-1}\bm{e}_{p-1} = \Vr_{p-1}$.
Our W-cycle algorithm terminates the recursion at level $k=1$
and directly solves $\bbA_{1}\bm{e}_{1} = \Vr_{1}$ or $\tilde{\bbA}_{1}\bm{e}_{1} = \Vr_{1}$. 
In the two-dimensional case, we employ MATLAB’s direct solver (e.g., backslash operator) for an exact solution, 
while for the three-dimensional case, the system is resolved via MUMPS package \cite{MUMPS}.

\subsection{Smoothing property}
\label{sec6.1}
In this test,
we verify the pre-smoothing property \eqref{pre-smoothing} of the polynomial smoother for linear system 
$\bbA_{k} \Vu_{k} = \Vf_{k} $.
The two-dimensional domain is given by $\Omega = (0,1)^{2} .$ 
The smoothing ratio is defined by 
\[
\Cs :=
\frac{\N{ \mathbb{M}^{-\frac{1}{2}}_k\bbA_k P_{m}\left(\omega_k \mathbb{D}_k^{-1}\bbA_k \right) \left(\Vu_{k} - \Vu_{k} ^{m} \right) }_{\bbA_k,0}}{\N{\Vu_{k} - \Vu_{k} ^{m} }_{\bbA_k,1} },
\]
obtained with $m$ times 
fourth-kind Chebyshev polynomial iterations, see \eqref{4th_Chebyshev}.

\begin{figure}[!htbp]
	\centering
	\includegraphics[width=4.5cm]{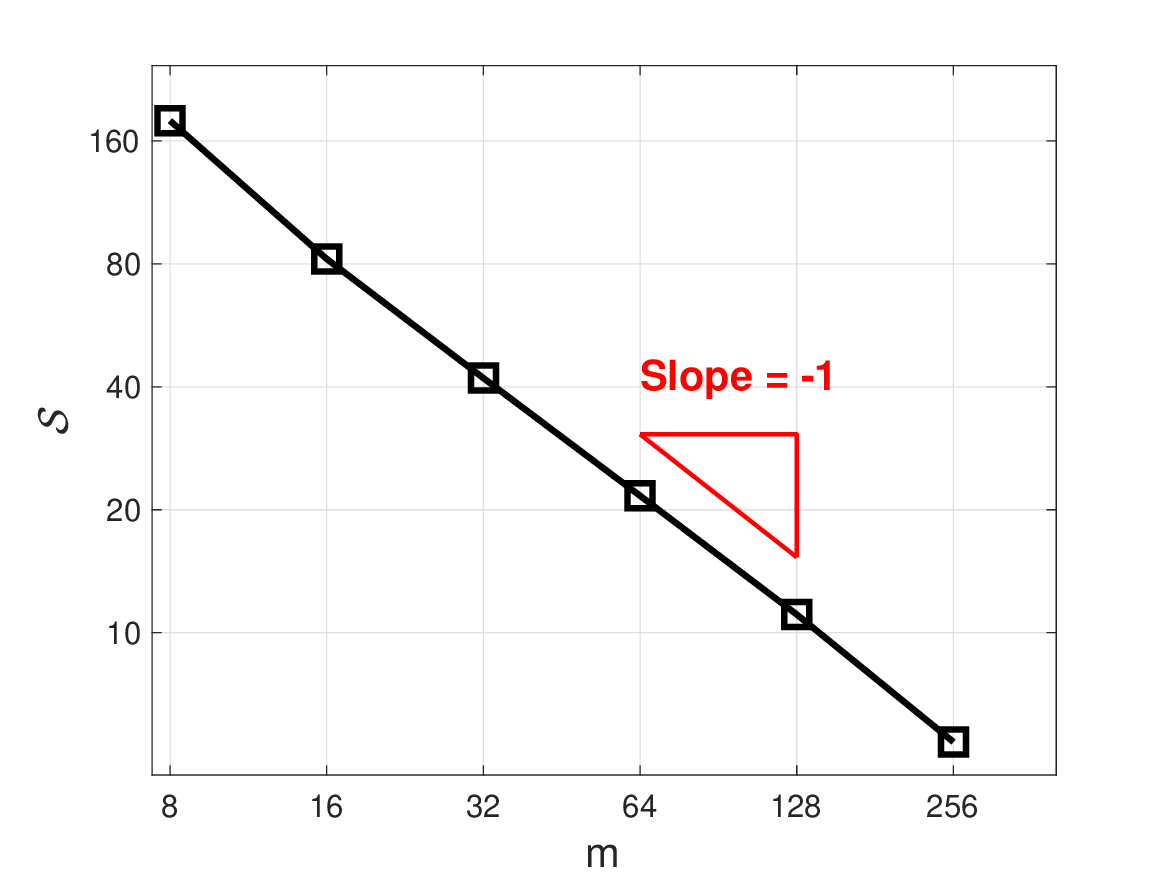}
	\includegraphics[width=4.5cm]{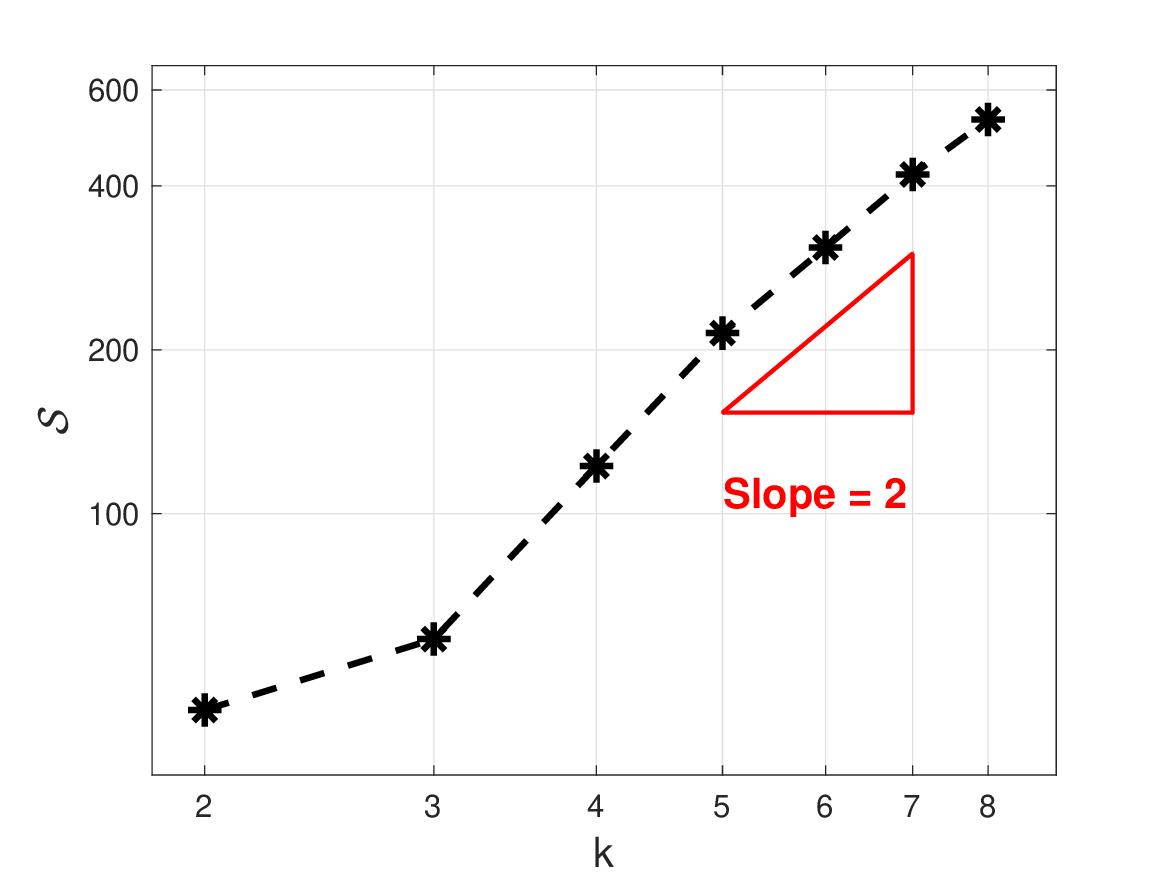}
	\includegraphics[width=4.5cm]{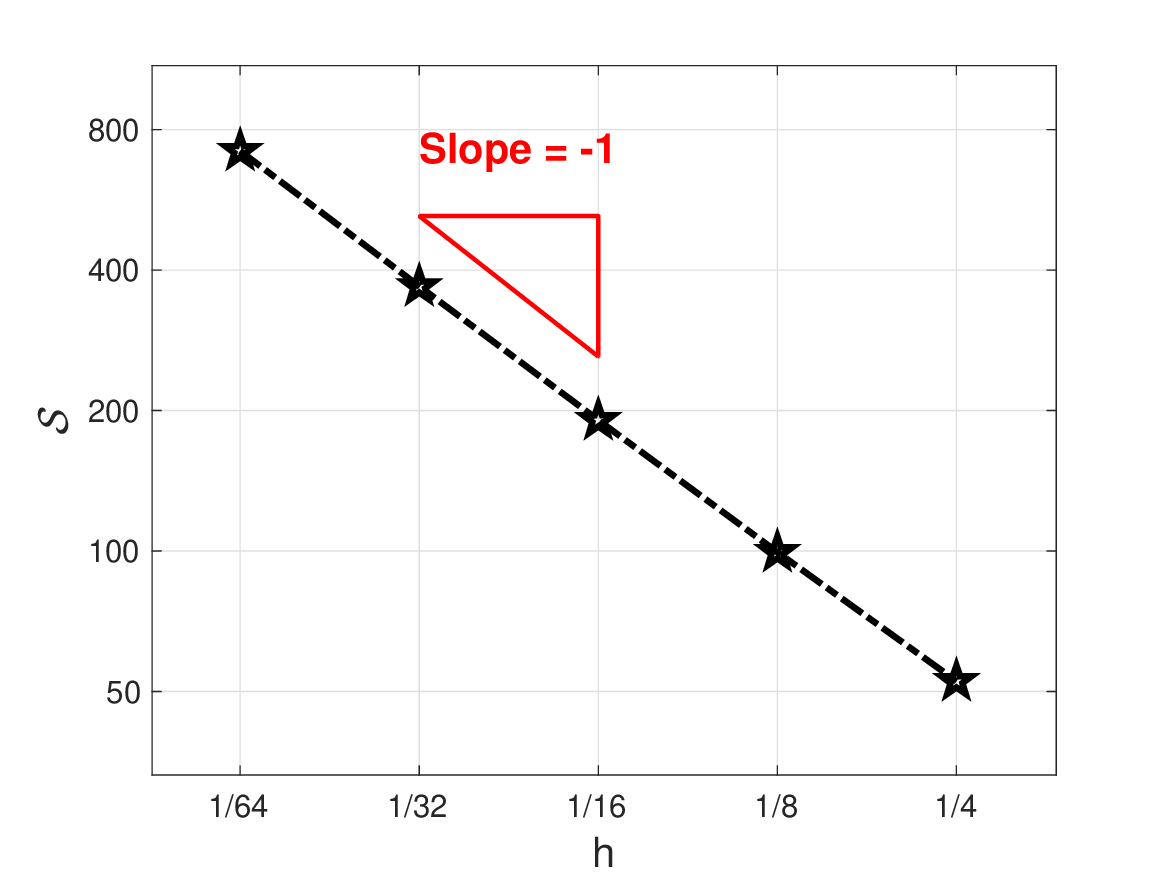}
	\caption{Estimates of the smoothing property in \eqref{pre-smoothing} with respect to $m$, $k$, and $h$.}
	\label{Fig_smoothing}
\end{figure}

Figure \ref{Fig_smoothing} illustrates the relationship of the smoothing ratio 
\( \Cs \) and  
smoothing times \( m \), degree of polynomial \( k \) and mesh size $h$.
The left subfigure shows that, with fixed \( k = 5 \) and \( N_h = 32^2 \) 
 \(( N_h = {1} / {h^2} )\), 
the smoothing ratio satisfies the asymptotic behavior \( \Cs \sim  m^{-1} \).
The middle subfigure demonstrates that, with fixed \( m = 20 \) and \( N_h = 32^2 \), 
the smoothing ratio satisfies the asymptotic behavior 
  $\Cs \sim  k^{2} $.
Finally, the right subfigure indicates that, with fixed \( k = 5 \) and \( m = k^{2}  \), 
the smoothing ratio satisfies the asymptotic behavior $\Cs \sim  h^{-1} $.
Together, these results confirm the theoretical prediction that $\Cs  \sim { k^{2}  }{ (hm)^{-1} }.$

\subsection{Two-level algorithm}
\label{sec6.2}
This example verifies the convergence rate of the two-level algorithm with non-inherited bilinear forms.
The domain is set by $\Omega = (-1,1) $.
We set the tolerance for relative residual of the two-level algorithm by $\varepsilon = 10^{-8}$,
which means that the iteration stops when
${||\Vf_p - \bbA_{p}\Vu^{N_{\rm it} }_h ||} / {||\Vf_p - \bbA_{p}\Vu^0_h ||} 
\leq \varepsilon$.
We define the convergence rate $\rho$ as 
\[
  \rho := \exp \left(\frac{1}{N_{\rm it} } 
   \ln \frac{||\Vf_p - \bbA_{p}\Vu^{N_{\rm it}}_h ||}{||\Vf_p - \bbA_{p}\Vu^0_h ||} \right).
\]

Table \ref{Tab_Acc1D_2level_poly} 
shows the convergence rate with different smoothing times 
$m = \left[p^{0.5}\right]$, $p$ and $\left[p^{1.5}\right]$.
It can be observed that for $ m = p $, 
the convergence rate is stable as $ p $ increases.
This is consistent with our theoretical results \eqref{Twoconv}.
Besides that, the numerical results show that the convergence rates are independent of mesh size $h$.

\begin{table}[!htbp]
\caption{Convergence rate $\rho$ of Two-level algorithm with different smoothing steps.} 
\centering
\begin{tabular}{c|c|c|c|c|c|c|c}
\toprule
\multicolumn{8}{c}{$m= \left[p^{0.5}\right]$} \\
\midrule
$N_h$ & {$p=2$} & {$p=3$} & {$p=4$} & {$p=5$} & {$p=6$} & {$p=7$} & {$p=8$}  \\ 
\hline
4  &  0.91  & 0.95 & 0.91 & 0.94 & 0.95 & 0.96  & 0.97 \\ \hline
8  &  0.91  & 0.95 & 0.92 & 0.94 & 0.95 & 0.96  & 0.97 \\ \hline
16 &  0.91  & 0.95 & 0.92 & 0.94 & 0.95 & 0.96  & 0.97 \\ \hline
32 &  0.91  & 0.95 & 0.92 & 0.94 & 0.95 & 0.96  & 0.97 \\ \hline
64 &  0.91  & 0.95 & 0.91 & 0.94 & 0.95 & 0.96  & 0.97 \\
\midrule
\multicolumn{8}{c}{$m = p$} \\
\midrule
$N_h$ & {$p=2$} & {$p=3$} & {$p=4$} & {$p=5$} & {$p=6$} & {$p=7$} & {$p=8$}  \\ 
\hline
4   & 0.79   & 0.80  & 0.75  & 0.75  & 0.72  & 0.72  & 0.71 \\ \hline
8   & 0.79   & 0.80  & 0.76  & 0.75  & 0.72  & 0.72  & 0.71 \\ \hline
16  & 0.79   & 0.80  & 0.76  & 0.75  & 0.72  & 0.72  & 0.71 \\ \hline
32  & 0.79   & 0.80  & 0.76  & 0.75  & 0.72  & 0.72  & 0.71 \\ \hline
64  & 0.79   & 0.80  & 0.76  & 0.75  & 0.72  & 0.72  & 0.71 \\
\midrule
\multicolumn{8}{c}{$m = \left[p^{1.5} \right]$} \\
\midrule
$N_h$ & {$p=2$} & {$p=3$} & {$p=4$} & {$p=5$} & {$p=6$} & {$p=7$} & {$p=8$}  \\ 
\hline
4   & 0.80  & 0.56  & 0.35  & 0.28  & 0.20  & 0.17  & 0.15 \\ \hline
8   & 0.80  & 0.58  & 0.36  & 0.28  & 0.20  & 0.17  & 0.15 \\ \hline
16  & 0.80  & 0.58  & 0.36  & 0.28  & 0.20  & 0.17  & 0.15 \\ \hline
32  & 0.80  & 0.58  & 0.36  & 0.28  & 0.20  & 0.17  & 0.14 \\ \hline
64  & 0.80  & 0.58  & 0.36  & 0.28  & 0.20  & 0.17  & 0.14 \\
\bottomrule
\end{tabular} 
\label{Tab_Acc1D_2level_poly}
\end{table}%


\subsection{W-cycle algorithm}
\label{sec6.3}
In this example, we consider the $W$-cycle algorithm of $p$-MG method.
We denote by $ K ( 2 \le K \leq p) $ the number of total levels,
$K=p$ implies that the coarsest level space uses linear broken polynomial space $U_{1}(\Ct_{h})$,
$ K = 2 $ indicates the two-level algorithm, 
and other cases mean that the coarsest level space corresponds to broken polynomial space $U_{p - K + 1}(\Ct_{h} )$.
The domain is set by $\Omega = (0,1)^{d}, d = 1,2,3$.

The iteration steps $N_{\rm it} $ 
and convergence rate $\rho$ for $K=p$ are list in Table \ref{Tab_Wcycle_poly}.
Notice that, at each level $k$, 
we take constant smoothing times $m=p$ for the inherited bilinear forms
and take changeable smoothing times $m=k$ for the non-inherited bilinear forms.
It can be observed that, regardless of whether inherited or non-inherited bilinear forms are used, 
the convergence rates are uniform with respect to $ p $ in dimensions $d\leq 3$, 
which also confirms the theoretical results provided in \eqref{W-cycle}. 
Since there is no big difference between the two approaches,  
we will only implement the inherited bilinear form in the following numerical tests.

\begin{table}[!htbp]
\centering
\caption{Convergence rate $\rho$ of W-cycle $p$-MG method with different bilinear forms $(K=p)$.}
\begin{tabular}{c|c|c|c|c|c}
\toprule
\multicolumn{6}{c}{1D Poisson problem} \\
\hline
{\textbf{bilinear forms}} & {$N_h$ } & {$p=3$} & {$p=4$} & {$p=5$} & {$p=6$} \\ \hline
\multirow{3}{*}{\textbf{inherited}} 
& 8  & 0.83 & 0.81 & 0.80 & 0.78 \\ \cline{2-6} 
& 16 & 0.83 & 0.81 & 0.80 & 0.78 \\ \cline{2-6} 
& 32 & 0.83 & 0.81 & 0.80 & 0.78 \\ \hline
\multirow{3}{*}{\textbf{non-inherited}} 
& 8  & 0.81 & 0.80 & 0.80 & 0.78 \\ \cline{2-6} 
& 16 & 0.81 & 0.80 & 0.80 & 0.78 \\ \cline{2-6} 
& 32 & 0.81 & 0.80 & 0.80 & 0.78 \\ 
\midrule
\multicolumn{6}{c}{2D Poisson problem} \\
\midrule
{\textbf{bilinear forms}} & {$N_h$ } & {$p=3$} & {$p=4$} & {$p=5$} & {$p=6$} \\ \hline
\multirow{3}{*}{\textbf{inherited}} 
& $8^2$    & 0.83   & 0.80   & 0.79   & 0.77 \\ \cline{2-6} 
& $16^2$   & 0.83   & 0.80   & 0.79   & 0.76 \\ \cline{2-6} 
& $32^2$   & 0.83   & 0.80   & 0.79   & 0.76 \\ \hline
\multirow{3}{*}{\textbf{non-inherited}} 
& $8^2$   & 0.84  & 0.81  & 0.79  & 0.78 \\ \cline{2-6} 
& $16^2$  & 0.83  & 0.80  & 0.79  & 0.78 \\ \cline{2-6} 
& $32^2$  & 0.83  & 0.80  & 0.78  & 0.78 \\
\midrule
\multicolumn{6}{c}{3D Poisson problem} \\
\midrule
{\textbf{bilinear forms}} & {$N_h$ } & {$p=3$} & {$p=4$} & {$p=5$} & {$p=6$} \\ \hline
\multirow{3}{*}{\textbf{inherited}} 
& $8^3$  & 0.75 & 0.72 & 0.70 & 0.68 \\ \cline{2-6} 
& $16^3$ & 0.75 & 0.72 & 0.70 & 0.68 \\ \cline{2-6} 
& $32^3$ & 0.75 & 0.72 & 0.70 & 0.68 \\ \hline
\multirow{3}{*}{\textbf{non-inherited}} 
& $8^3$    & 0.76   & 0.73 & 0.71  & 0.68 \\ \cline{2-6} 
& $16^3$   & 0.76   & 0.73   & 0.71  & 0.68 \\ \cline{2-6} 
& $32^3$   & 0.76   & 0.73   & 0.71  & 0.68 \\ 
\bottomrule
\end{tabular}
\label{Tab_Wcycle_poly}
  \end{table}

\subsection{Operator complexity}
\label{sec6.4}
In this example, we take $K=p$ and investigate operator complexity $C_O$ to evaluate the computational overhead of a multigrid method \cite{Stuben1999},
\[
C_O := {\sum_{k=1} ^{p}  nnz(\bbA_{k} )}/{nnz(\bbA_p) }.
\] 
Here $nnz(\bbA_{k}) $ denotes the number of nonzero entries in $\bbA_k$.

The results of the operator complexity with different dimensions are shown in Table \ref{Tab_OC}.
It shows that our $p$-MG method maintains stable operator complexity with respect to the mesh size, 
highlighting its scalability and efficiency in dealing with large-scale problems. 
This means the computational complexity not increase significantly with mesh refinements.  
The operator complexity does increase with large $ p $, 
this is expected due to complex polynomial spaces of high orders, 
leading to additional computations. 
However, even with increasing $ p $, the growth in complexity remains moderate, 
demonstrating the robustness and efficiency of the algorithm dealing with high-order methods. 

\begin{table}[!htbp]
\centering
\caption{Operator complexity $C_{O} $ of the $p$-MG method $(K=p)$.}
\begin{tabular}{c|c|c|c|c|c}
\toprule
\multicolumn{6}{c}{1D Poisson problem} \\ 
\midrule
$N_h$   & $p=2$ & $p=3$ & $p=4$ & $p=5$ & $p=6$ \\ \hline
4 & 1.44 & 1.82 & 2.16 & 2.51 & 2.84 \\ \hline
8 & 1.44 & 1.82 & 2.16 & 2.51 & 2.84 \\ \hline
16 & 1.44 & 1.82 & 2.16 & 2.51 & 2.84 \\ 
\midrule
\multicolumn{6}{c}{2D Poisson problem} \\ 
\midrule
$N_h$ & $p=2$ & $p=3$ & $p=4$ & $p=5$ & $p=6$ \\ \hline
$4^2$ & 1.28 & 1.54 & 1.77 & 2.02 & 2.26 \\ \hline
$8^2$ & 1.28 & 1.54 & 1.77 & 2.02 & 2.27 \\ \hline
$16^2$ & 1.28 & 1.54 & 1.77 & 2.02 & 2.27 \\ 
\midrule
\multicolumn{6}{c}{3D Poisson problem} \\ 
\midrule
$N_h$ & $p=2$ & $p=3$ & $p=4$ & $p=5$ & $p=6$ \\ \hline
$4^3$ & 1.06 & 1.20 & 1.30 & 1.45 & 1.54 \\ \hline
$8^3$ & 1.07 & 1.20 & 1.30 & 1.45 & 1.54 \\ \hline
$16^3$ & 1.07 & 1.20 & 1.30 & 1.43 & 1.54 \\ 
\bottomrule
\end{tabular}
\label{Tab_OC}
\end{table}

\subsection{Convergence rates with respect to $m$ and $p$}
\label{sec6.5}
In this example, we compare our $p$-MG method with that in \cite{Antonietti2015}.
We fix the mesh size $h=0.0625$, penalty parameter $\alpha_0=10$ in the two-dimensional case, 
and analyze the relationship between convergence rates and $ p $, $m$ and $ K $.  

First, we fix $p = 5$ and compare the convergence rates for different smoothing times and grid levels. 
As shown in Table \ref{Tab_rho_vs_m_K}, 
the convergence rate remains stable to the number of levels $K$,
and decreases significantly $m$ increases.
The ``Ref'' columns in Table indicate that the corresponding datum are obtained \cite{Antonietti2015}, 
while the ``our'' columns stand for our numerical results.

\begin{table}[!htbp]
\centering
\caption{Convergence rate of the $p$-MG method with respect to $m$ and $K$ $(p = 5, h=0.0625)$.}
\begin{tabular}{c|c|c|c|c|c|c|c|c|c|c}
\toprule
\multirow{2}{*}{$K$}
& \multicolumn{2}{c|}{$m=2$} & \multicolumn{2}{c|}{$m=4$} & \multicolumn{2}{c|}{$m=6$} & \multicolumn{2}{c|}{$m=10$} & \multicolumn{2}{c}{$m=20$} \\ \cline{2-11}
& Ref  & our   & Ref   & our   & Ref   & our   & Ref   & our   & Ref   & our    \\ \hline
2 & 0.93 & 0.94  & 0.88 & 0.85  & 0.85 & 0.71  & 0.82 & {\bf 0.41}  & 0.73 & {\bf 0.06}  \\ 
\hline
3  & 0.93 & 0.94  & 0.88 & 0.85  & 0.85 & 0.71  & 0.82 & {\bf 0.41}  & 0.73 & {\bf 0.06}  \\
\hline
4  & 0.95 & 0.96  & 0.92 & 0.85  & 0.89 & 0.71  & 0.84 & {\bf 0.41}  & 0.73 & {\bf 0.06}  \\
\bottomrule
\end{tabular}
\label{Tab_rho_vs_m_K}
\end{table}

Next, we fix $m=10$ and enlarge $ p $ from $2$ to $6$.
It can be observed from Table \ref{Tab_rho_vs_k_K} that the convergence rate remains 
stable to the number of total levels $K$.
As $ p $ increases, the number of degrees of freedom increase significantly,
and lead to an increase of the convergence rate.

\begin{table}[!htbp]
\centering
\caption{Convergence rate of the $p$-MG method with respect to $p$ and $K$ $(m = 10, h=0.0625)$.}
\begin{tabular}{c|c|c|c|c|c|c|c|c|c|c}
\toprule
\multirow{2}{*}{$K$ }
& \multicolumn{2}{c|}{$p=2$} & \multicolumn{2}{c|}{$p=3$} & \multicolumn{2}{c|}{$p=4$} & \multicolumn{2}{c|}{$p=5$} & \multicolumn{2}{c}{$p=6$} \\ \cline{2-11}
& Ref  & our   & Ref   & our   & Ref   & our   & Ref   & our   & Ref   & our    \\ \hline
2  & 0.58 & {\bf 0.06}  & 0.79 & {\bf 0.20}  & 0.77 & {\bf 0.28}  & 0.82 & {\bf 0.41}  & 0.85 & {\bf 0.49}  \\
\hline
3  & -    & -  		& 0.78  & {\bf 0.20} & 0.78  & {\bf 0.28} & 0.82  & {\bf 0.41} & 0.85  & {\bf 0.49}         \\
\hline
4  & -       & -      		& -       & - 		   &  0.78  & {\bf 0.29} & 0.84  & {\bf 0.41} & 0.84  & {\bf 0.49}         \\
\bottomrule
\end{tabular}
\label{Tab_rho_vs_k_K}
\end{table}

Furthermore, we employ the inverse of damped diagonal matrix of matrix $\bbA_{k} $, 
defined as \( \omega_k \bbD^{-1}_k \), in smoothing steps.
Since computing the inverse of \( \bbD_k \) is very cheap,
it does not yield additional computational cost. 
Consequently, the computational overhead for our \( p \)-multigrid method 
remains essentially the same in \cite{Antonietti2015}. 
Overall, the numerical results demonstrate that our \( p \)-multigrid method is efficient, 
particularly for not very small $m$, yields a significantly fast convergence rate.

\subsection{$p$-MG preconditioned GMRES method}
\label{sec6.6}
In this example, 
We test the $p$-MG preconditioned GMRES method for $p=2,3,4,5,6$ on different meshes.
We employ inherited bilinear form W-cycle $p$-MG as the preconditioner, 
where only a single iteration is performed at each preconditioning step.
The relative tolerance of GMRES method is $10^{-8}$.
Other settings for this example is same as example \ref{sec6.3}.
$N_{\rm GMRES} $ denotes the number of iteration steps 
for the preconditioned GMRES iterative methods. 

First, we uniformly refine an uniform mesh.
As demonstrated in Table \ref{pmg_GMRES}, 
$p$-MG preconditioned GMRES method is stable independent of the mesh size and polynomial degree.
This stands in sharp contrast to the standard GMRES method (without preconditioning).

\begin{table}[!htbp]
\centering
\caption{$N_{\rm GMRES}$ of $p$-MG preconditioned GMRES methods on uniform mesh.} 
\begin{tabular}{c|c|c|c|c|c||c|c|c|c|c}
\toprule
\multirow{2}{*}{ $N_h$ }
& \multicolumn{5}{c||}{GMRES} & 
\multicolumn{5}{c}{ $p$-MG preconditioned GMRES } \\ \cline{2-11}
& $p=2$ & $p=3$ & $p=4$ & $p=5$ & $p=6$ & $p=2$ & $p=3$ & $p=4$ & $p=5$ & $p=6$\\ \hline
$8^{3} $   & 222 & 335  &537  & 668 & 936 	& \textbf{17 }& \textbf{15} & \textbf{15} & \textbf{14} & \textbf{13} \\ \hline
$16^{3} $  & 404 & 599  &921  & 1089& 1396 	& \textbf{17} & \textbf{16} & \textbf{15} & \textbf{14} & \textbf{13}\\ \hline 
$32^{3} $  & 764 & 1041 &1902 & 2097& 3225  & \textbf{18} & \textbf{15} & \textbf{15} & \textbf{14} & \textbf{13} \\ 
\bottomrule
\end{tabular}
\label{pmg_GMRES}
\end{table}

Next, we uniformly refine two adapted meshes.
As demonstrated in Table \ref{pmg_GMRES_2}, 
$p$-MG preconditioned GMRES method is stable independent of the mesh size and polynomial degree.

\begin{table}[!htbp]
\centering
\caption{$N_{\rm GMRES}$ of $p$-MG preconditioned GMRES method on adopted mesh.} 
\begin{tabular}{c|c|c|c|c|c||c|c|c|c|c|c}
\toprule
{ $N_h$ }
& $p=2$ & $p=3$ & $p=4$ & $p=5$ & $p=6$ & $N_{h} $ & $p=2$ & $p=3$ & $p=4$ & $p=5$ & $p=6$\\ \hline
15 	& 15 & 16 & 16 & 15 & 14 & 22		& 15 & 16 & 16 & 15 & 14 \\ \hline
120 	& 17 & 17 & 16 & 15 & 14 & 176	 	& 17 & 17 & 16 & 15 & 14 \\ \hline
960 	& 18 & 17 & 16 & 15 & 14 & 1408 	& 18 & 17 & 16 & 15 & 14 \\
\bottomrule
\end{tabular}
\label{pmg_GMRES_2}
\end{table}

\section{Elliptic interface problems}
\label{sec7}
In this section, we investigate $p$-MG method for unfitted finite element method solving elliptic interface problems. 

Interface problems occur in various physical and engineering applications. 
Body-fitted finite element methods effectively address the geometry of these interfaces by ensuring that the vertices of the finite element mesh align precisely along them 
\cite{Babuska70,Chen98}. 
However, in domains with complex shapes, constructing body-fitted shape-regular finite element meshes can be both challenging and time-consuming. 
Unfitted finite element methods within the framework of discontinuous Galerkin (DG) have attracted considerable attention in scholarly discourse over the past twenty years \cite{Hansbo02}.

A notable challenge in the application of unfitted finite element methods is the so-called small cut cell problem: 
these cut cells can be extraordinarily small and anisotropic, 
which may lead to severe ill-conditioning of the stiffness matrix, 
particularly in the context of high-order finite element methods.
Merging the small cut elements with neighboring large elements is appropriate technique to address this issue \cite{Chen2021, Chen2023, Chen2024}.

Let $\Omega=\bar\Omega_{1} \cup \bar\Omega_{2}, \Omega_{1} \cap \Omega_{2} = \emptyset, \Gamma = \partial\Omega_{1} \cap \partial \Omega_{2} $. 
This algorithm focuses on the merging elements into macro elements and the strategically subdividing elements to ensure that the resulting mesh complies with the Big Element Conditions (BEC):
\begin{itemize}
\item At most one cut on each element edge and no face-only cuts: $h \kappa \le \textsf{ctol}$ ($\kappa$ := curvature);
\item Interface deviation control: $\eta_{T} \le \textsf{dtol} $\;($T:=$ element, $\eta_{T} $ := interface deviation in \cite{Chen2023});
\item Edge-cut ratio control: $|e \cap \Omega| / |e| \ge \textsf{etol}\;(e := \text{edge})$;
\item (optional) Volume ratio control: $|T \cap \Omega_{i} | / |T| \ge \textsf{vtol}, i = 1,2$;
\end{itemize}
where the induced mesh is composed macro and unmerged elements.
For the 2D case, a provable algorithm for reliably merging rectangle elements into 
rectangle macro elements satisfying the BEC, for arbitrary locally refine non-conforming rectangular meshes, has been proposed in \cite{Chen2023}.

PHG \cite{Zhang2009} has implemented an element merging algorithm which can generate an induced mesh fully satisfying BEC on any locally adapted cuboidal mesh.
First, we set $\Omega=[0,1]^3$, and test the following three interfaces $\Gamma_i,\, i=1,2,3$.
The three interfaces are a torus, an ellipsoid, and a sphere, respectively.
\begin{align*}
  \Gamma_1 & = \left\{ (\sqrt{(x-0.5)^2+(y-0.5)^2} -0.25)^2 + (z-0.5)^2 = 0.1^2 \right\}, \\
  \Gamma_2 & = \left\{(x-0.5)^2+(y-0.5)^2+(z-0.5)^2=0.3^2\right\}, \\
  \Gamma_3 & = \left\{\frac{ (x-0.5)^2 }{ 0.2^2 } + \frac{ (y-0.5)^2 }{ 0.4^2 } + \frac{ (z-0.5)^2 }{ 0.3^2 } = 1\right \}.
\end{align*}

\begin{figure}[!htbp]
	\centering
	\includegraphics[width=3.5cm]{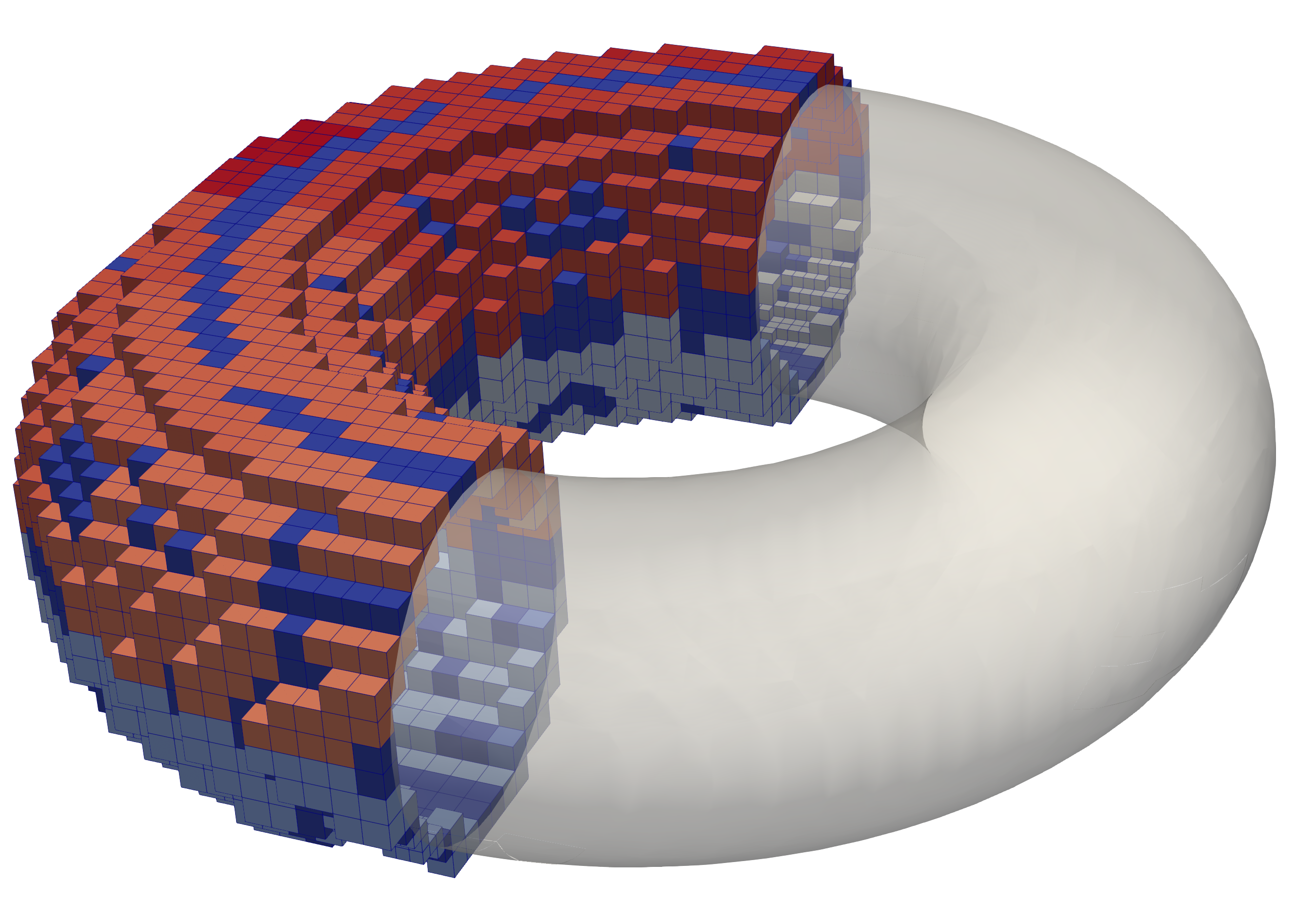}
  \hspace{15pt}
	\includegraphics[width=2.8cm]{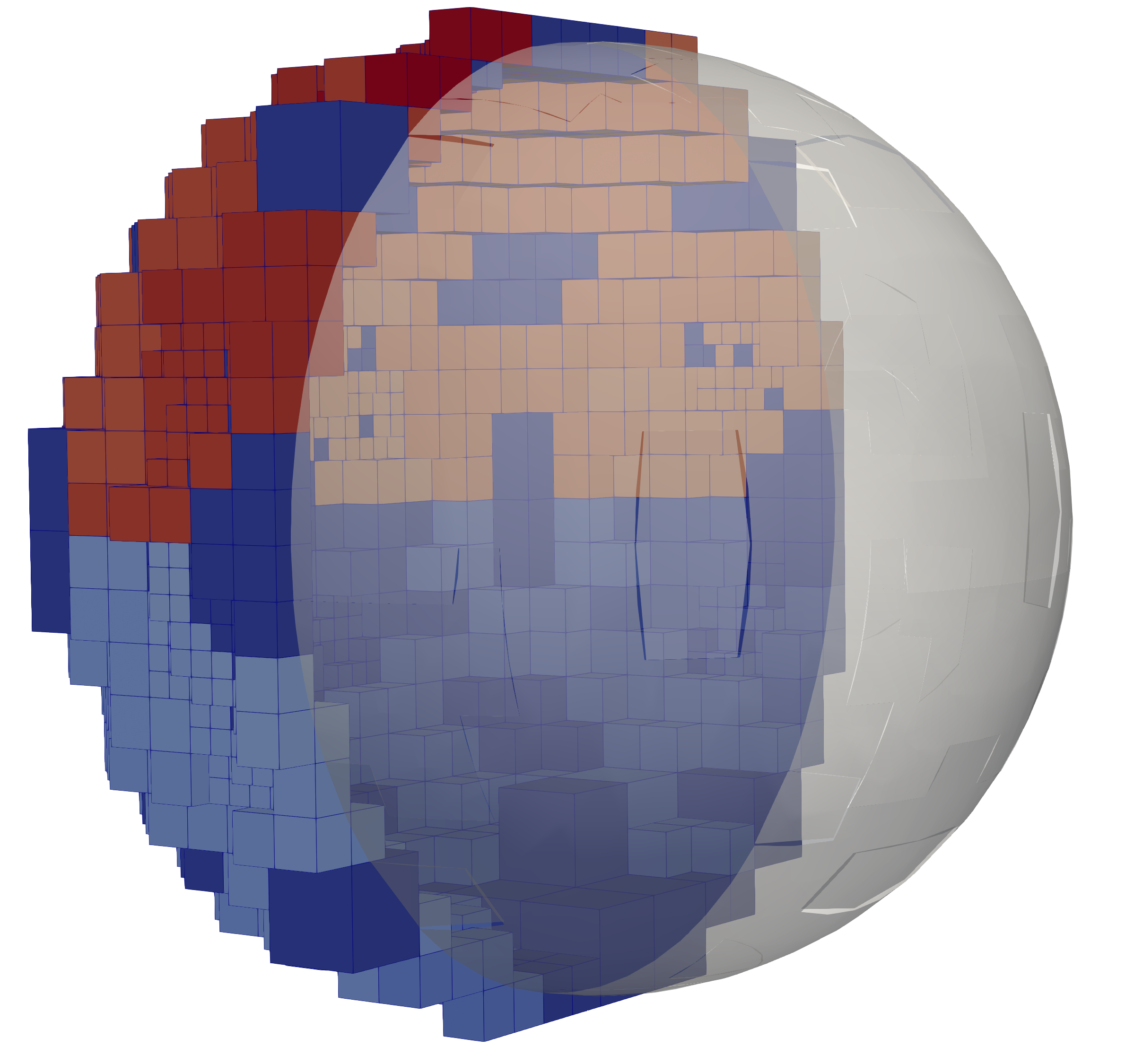}
  \hspace{15pt}
	\includegraphics[width=3cm]{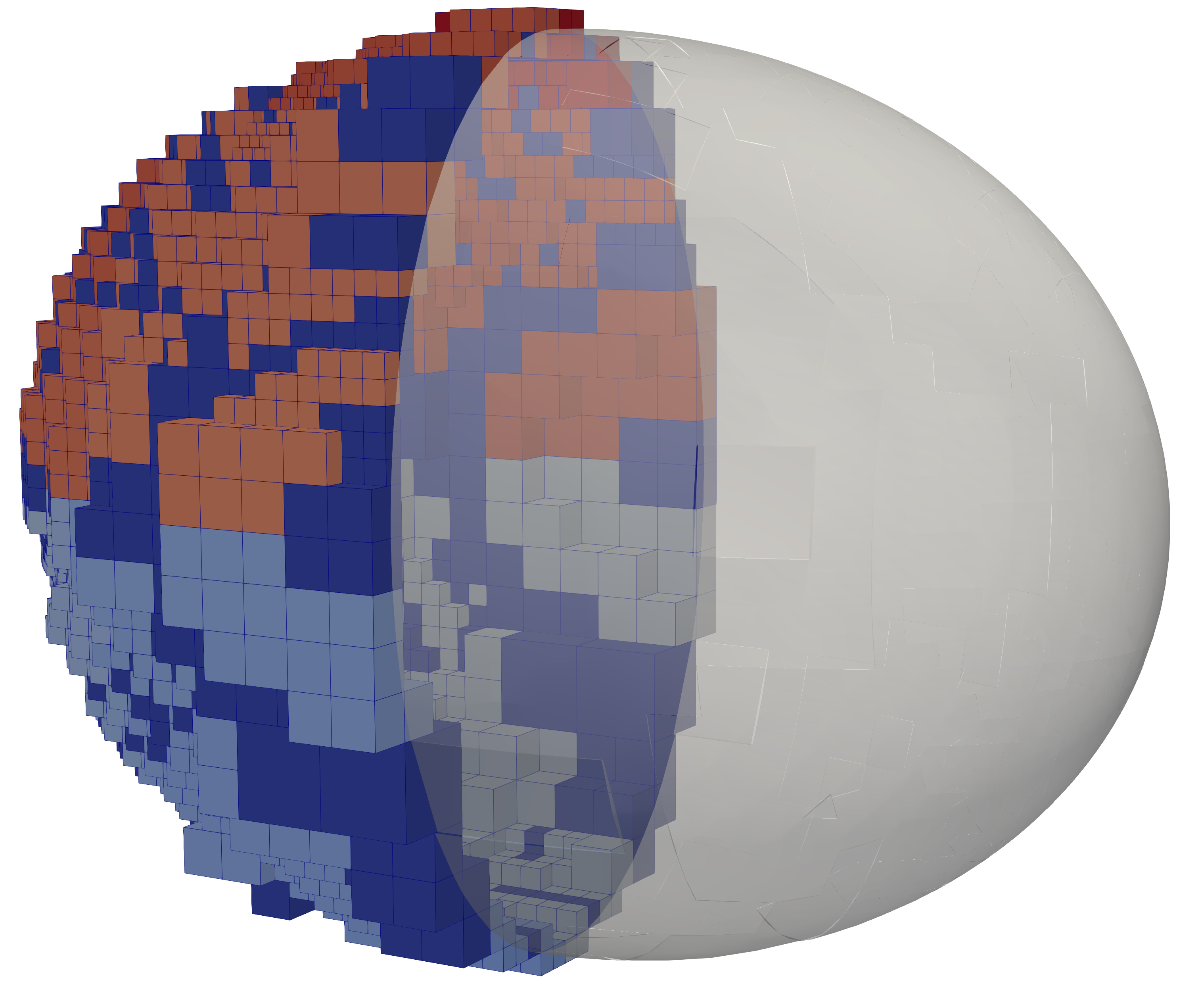}
	\caption{Macro elements by merging algorithm for different interfaces $\Gamma_{i}, 
	i=1,2,3$ with $r=3, p=2$.}
\label{Fig_macors}
\end{figure}

\begin{figure}[!htbp]
	\centering
	\includegraphics[width=3cm]{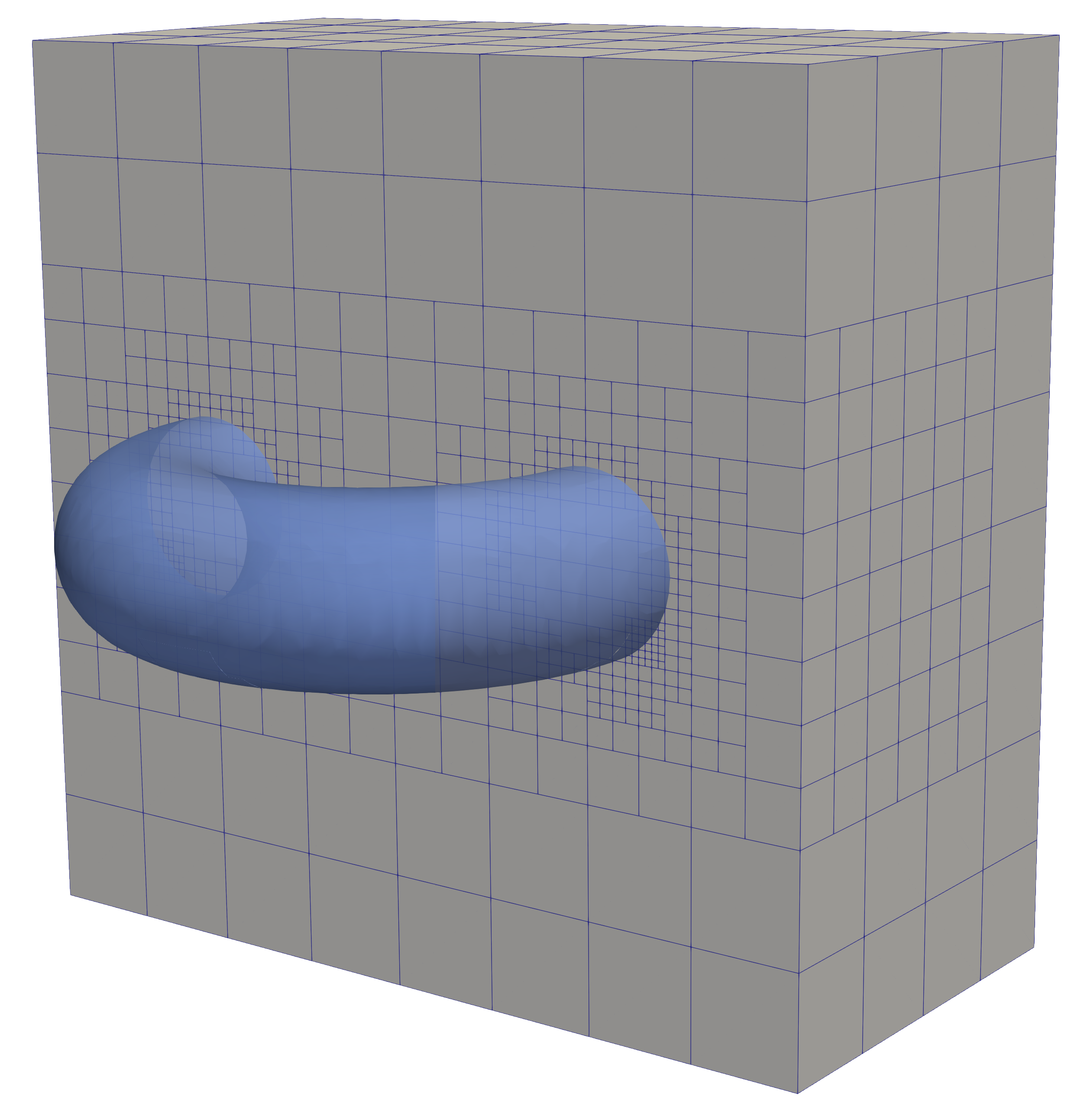}
  \hspace{15pt}
	\includegraphics[width=3cm]{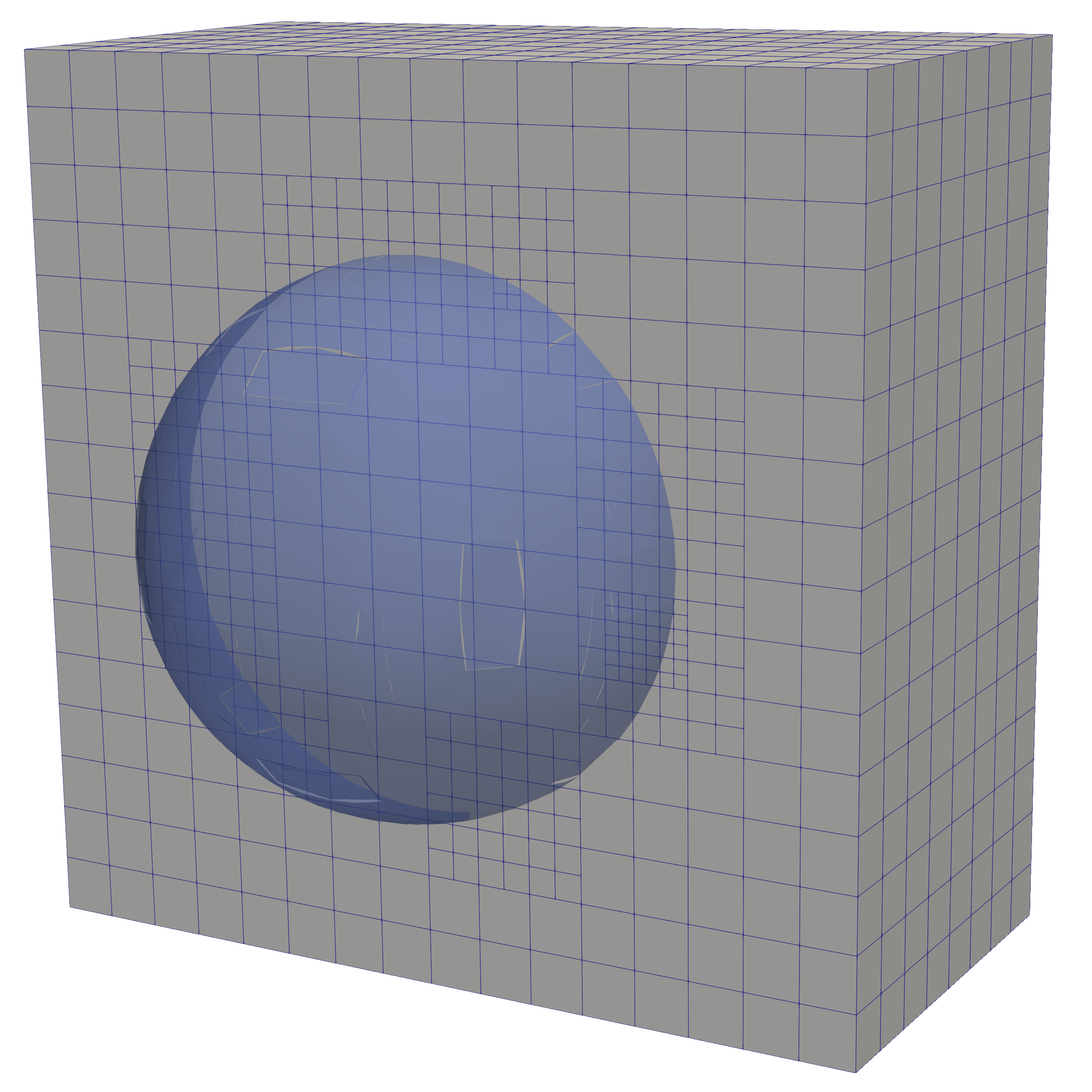}
  \hspace{15pt}
	\includegraphics[width=3cm]{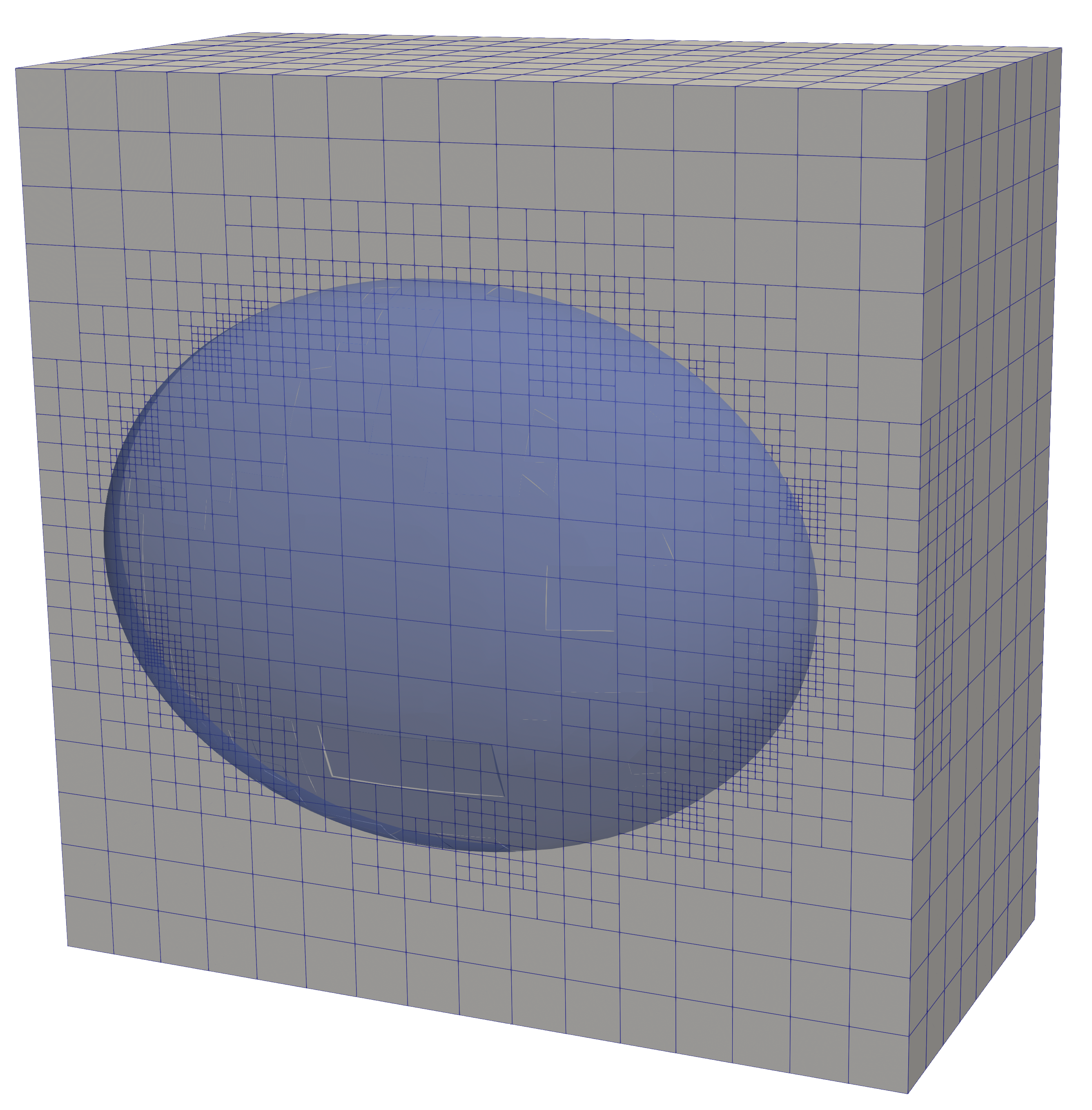}
	\caption{Induced meshes by merging algorithm for different interfaces $\Gamma_{i}, 
	i=1,2,3$, with $r=3, p=2$.}
	\label{Fig_induced_mesh}
\end{figure}

We use a cube as initial coboidal mesh denoted by $\Cm_{0}$, and take $r$ times uniform refine to generate new mesh denoted by $\Cm_{r}$. 
Then we use merging algorithm on mesh $\Cm_{r} $ for different interfaces.
The parameters of BEC for merging algorithm are given in Table \ref{Tab_Para_BEC}. 
Figure \ref{Fig_macors} and \ref{Fig_induced_mesh} show the macro elements and induced meshes by the merging algorithm for different interfaces with $r = 3, p = 2$.

\begin{table}[!htbp]
\centering
\caption{Parameters of BEC for merging algorithm.}
\begin{tabular}{c|c|c|c|c|c}
\toprule
\multirow{2}{*}{\textsf{ctol}} & \multicolumn{3}{c|}{\textsf{dtol}} & \multirow{2}{*}{\textsf{etol}} & \multirow{2}{*}{\textsf{vtol}} \\ \cline{2-4}
& $\Gamma_{1} $& $\Gamma_{2} $ & $\Gamma_3$& &  \\ \hline
0.5 & 0.2 & 0.1 & 0.08 &  0.25 & 0.0 \\
\bottomrule
\end{tabular}
\label{Tab_Para_BEC}
\end{table}

The DG formulation of the elliptic interface problem is given in \cite{Chen2021} and omitted here.
We test the $p$-MG method preconditioned GMRES method for different $p$, and $r$. 
Other setting is same as example \ref{sec6.6}.
Table \ref{Tab_pmg_gmres_interface} shows that
$p$-MG preconditioned GMRES method is stable with respect to $p$ and $r$ on the same interface.

\begin{table}[!htbp]
\centering
\caption{$N_{\rm GMRES} $ of $p$-MG preconditioned GMRES method for the DG scheme solving elliptic interface problems.}
\begin{tabular}{c|c|c|c|c|c|c}
\toprule
{Interfaces} & {$r$} & $N_{h}$ &  {$p=2$ } & {$p=3$ } & {$p=4$ } & {$p=5$ } \\ \hline 
\multirow{2}{*}{Torus }
&3 &25467& 30& 29 & 27 & 34 \\ \cline{2-7} 
&5 &50380& 29& 27 & 29 & 31 \\ \hline
\multirow{2}{*}{Sphere }
&3 &11208& 40  & 31  & 26 & 26\\ \cline{2-7}
&5 &37416& 40  & 29  & 24 &  23\\ \hline
\multirow{2}{*}{Ellipsoid }
&3 &22212& 58  & 70  & 47  & 45\\  \cline{2-7}
&5 &42316& 58  & 55  & 35 & 35  \\  
\bottomrule
\end{tabular}
\label{Tab_pmg_gmres_interface}
\end{table}

\section{Conclusion}
\label{sec8}
In this paper, we have developed an efficient $p$-MG method 
for high-order SIPDG discretization of the Poisson problem,
by using polynomial smoother based on the fourth-kind Chebyshev polynomial iterative method.
We provided a rigorous convergence analysis for the two-level and $W$-cycle algorithm of $p$-MG method, 
with non-inherited and inherited bilinear forms.
The convergence rate is independent of the mesh size $h$ and 
is uniform to polynomial degree $p$ when number of smoothing times satisfies $m \gtrsim p$.
Numerical results verify our theoretical findings, and also show that
$p$-MG method is stable for DG discretizations on locally refine meshes or unfitted meshes.
In future works, we plan to extend the study
to heterogeneous diffusion problems with discontinuous coefficients.
Additionally, we will also explore $\Hcurl$- and $\Hdiv$-elliptic (interfaces) problems 
and establish a rigorous convergence theory for these cases.

While the present work focuses on Cartesian grids, 
the proposed methodology can be naturally extended to triangular or tetrahedral meshes. 
However, such an extension necessitates a re-evaluation of the stiffness matrix's spectral properties, 
particularly its condition number, which introduces additional computational 
and theoretical challenges. A rigorous analysis of these aspects 
will be the subject of future research.

\bibliographystyle{plain}

\begin{thebibliography}{100}
\bibitem{MUMPS}
{\sc P. R. Amestoy, A. Buttari, J.-Y. L’Excellent and T. Mary},
{\it Performance and scalability of the block low-rank multifrontal
factorization on multicore architectures},
ACM Transactions on Mathematical Software,
45 (2019), pp 2:1--2:26.


\bibitem{Antonietti2020}
{\sc P. F. Antonietti and L. Melas},
{\it Algebraic multigrid schemes for high-order nodal discontinuous Galerkin methods},
SIAM J. Sci. Comput., 42 (2020), pp. A1147--A1173.

\bibitem{Antonietti2019}
{\sc P. F. Antonietti, G. Pennesi}, 
{\it V-cycle multigrid algorithms for discontinuous Galerkin methods on non-nested polytopic meshes}, 
J. Sci. Comput., 78 (2019), pp. 625--652.


\bibitem{Antonietti2017}
{\sc P. F. Antonietti, P. Houston, X. Hu, M. Sarti and M. Verani },
{\it Multigrid algorithms for $p$-version interior penalty discontinuous Galerkin methods on polygonal and polyhedral meshes}, 
Calcolo, 54 (2017), pp. 1169--1198. 

\bibitem{Antonietti2017_JSC}
{\sc P. F. Antonietti, M. Sarti, M. Verani, L. T. Zikatanov}, 
{\it A uniform additive Schwarz preconditioner for high-order discontinuous Galerkin approximations of elliptic problems}, 
J. Sci. Comput., 70 (2017), pp. 608--630.

\bibitem{Antonietti2015}
{\sc P. F. Antonietti, M. Sarti and M. Verani},
{\it Multigrid algorithms for hp-discontinuous Galerkin discretizations of elliptic problems},
SIAM J. Numer. Anal., 53 (2015), pp. 598--618.


\bibitem{Antonietti2011}
{\sc P. F. Antonietti, P. Houston},
{\it A class of domain decomposition preconditioners for hp-discontinuous Galerkin finite element methods},
J. Sci. Comput., 46 (2011), pp. 124--149. 

\bibitem{Arnold1982}
{\sc D. Arnold}, 
{\it An interior penalty finite element method with discontinuous elements},
SIAM J. Numer. Anal., 19 (1982), pp. 742--760.

\bibitem{Arnold2002}
{\sc D. Arnold, F. Brezzi, B. Cockburn and L. D. Marini},
{\it Unified analysis of discontinuous Galerkin methods for elliptic problems},
SIAM J. Numer. Anal., 02 (2002), pp. 1749--1779.

\bibitem{Babuska70}
{\sc I.~Babu\"{s}ka}, 
{\it The finite element method for elliptic equations with discontinuous coefficients}, {Computing}, 5 (1970), 207-213.

\bibitem{Bastian2012}
{\sc P. Bastian, M. Blatt and R. Scheichl},
{\it Algebraic multigrid for discontinuous Galerkin discretizations of heterogeneous elliptic problems},
Numer. Linear Algebra Appl., 19 (2012), pp. 367--388.

\bibitem{Brenner2005}
{\sc S. C. Brenner and J. Zhao},
{\it Convergence of multigrid algorithms for interior penalty methods},
Applied Numerical Analysis \& Computational Mathematics, 2 (2005), pp. 3--18.

\bibitem{Cangiani2014}
{\sc A. Cangiani, E. H. Georgoulis, P. Houston}, 
{\it hp-version discontinuous Galerkin methods on polygonal and polyhedral meshes}, 
Math. Models Methods Appl. Sci., 24 (2014), pp. 2009--2041.

\bibitem{Chen98}
{\sc Z.~Chen, J.~Zou}, 
{\it Finite element methods and their convergence for elliptic and parabolic interface problems}, {Numer. Math.}, 79, (2008), 175-202.

\bibitem{Chen2021}
{\sc Z. Chen, K. Li, X. Xiang},
{\it An adaptive high-order unfitted finite element method for elliptic interface problems}, 
Numer. Math., 149 (2021), pp. 507--548.

\bibitem{Chen2023}
{\sc Z. Chen, Y. Liu},
{\it An arbitrarily high order unfitted finite element method for elliptic interface problems 
with automatic mesh generation}, 
J. Comput. Phys., 491 (2023), 112384.

\bibitem{Chen2024}
{\sc Z.~Chen, K.~Li, M.~Lyu, and~X.~Xiang}, 
{\it A High Order Unfitted Finite Element Method for Time-Harmonic Maxwell Interface Problems}, 
{ Int. J. Numer. Anal. Model.}, 21 (2024), 822-849.

\bibitem{Dios2017}
{\sc B. Ayuso de Dios, R. Hiptmair and C. Pagliantini}, 
{\it Auxiliary space preconditioners for SIP-DG discretizations of H(curl)-elliptic problems with discontinuous coefficients}, 
IMA J. Numer. Anal., 37 (2017), pp. 646--686.

\bibitem{Dios2014}
{\sc B. Ayuso de Dios, M. Holst, Y. Zhu and L. Zikatanov}, 
{\it Multilevel preconditioners for discontinuous Galerkin
approximations of elliptic problems with jump coefficients}, 
Math. Comp., 83 (2014), pp. 1083--1120.

\bibitem{Dios2009}
{\sc B. Ayuso de Dios and L. Zikatanov},
{\it Uniformly convergent iterative methods for discontinuous Galerkin discretizations}, 
J. Sci. Comput., 40 (2009), pp. 4--36.

\bibitem{Epshteyn2007}
{\sc Y. Epshteyn and B. Rivière},
{\it Estimation of penalty parameters for symmetric interior penalty Galerkin methods},
J. Comput. Appl. Math., 206 (2007), pp. 843--872.

\bibitem{Gunatilake2022}
{\sc J. Gunatilake},
{\it An improved multigrid solver for the $p$-hierarchical basis finite element method using a space decomposition smoother},
Comput. Math. Appl., 124 (2022), pp. 52--62.

\bibitem{Hansbo02}
{\sc A.~Hansbo, and~P.~Hansbo}, 
{\it An unfitted finite element method based on Nitsche’s method for elliptic interface problems},{Comput. Methods Appl. Mech. Eng. }, 191 (2002), 5537-5552.

\bibitem{Hartmann2009}
{\sc R. Hartmann, M. Lukacova-Medvidova and F. Prill},
{\it Efficient preconditioning for the discontinuous Galerkin finite element method by low-order elements},
Appl. Numer. Math., 59 (2009), pp. 1737--1753.


\bibitem{Haupt2013}
{\sc L. Haupt, J. Stiller and W. E. Nagel},
{\it A fast spectral element solver combining static condensation and multigrid techniques},
J. Comput. Phys., 255 (2013), pp. 384--395.

\bibitem{Helenbrook2006}
{\sc B. T. Helenbrook and H. L. Atkins}, 
{\it Application of $p$-multigrid to discontinuous Galerkin formulations of the Poisson equation}, 
AIAA journal, 44 (2006), pp. 566--575.

\bibitem{Helenbrook2003}
{\sc B. T. Helenbrook, D. J. Mavriplis and H. L. Atkins}. 
{\it Analysis of ``p''-multigrid for continuous and discontinuous finite element discretizations}, 
16th AIAA Computational Fluid Dynamics Conference, (2003).


\bibitem{Hong2016}
{\sc Q. Hong, J. Kraus, J. Xu and L. Zikatanov},
{\it A robust multigrid method for discontinuous Galerkin discretizations of Stokes and linear elasticity equations},
Numer. Math., 132 (2016), pp. 23--49.

\bibitem{Lazarov2007}
{\sc R. D. Lazarov and S. D. Margenov}, 
{\it CBS constants for multilevel splitting of graph-Laplacian and application to preconditioning of discontinuous Galerkin systems}, 
J. Complexity., 23 (2007), pp. 498--515.

\bibitem{Lazarov2006}
{\sc V. A. Dobrev, R. D. Lazarov, P. S. Vassilevski and L. T. Zikatanov},
{\it Two-level preconditioning of discontinuous Galerkin approximations of second-order elliptic equations},
Numer. Linear Algebra Appl., 13 (2006), pp. 753--770.


\bibitem{Lottes2023}
{\sc J. Lottes},
{\it Optimal polynomial smoothers for multigrid V-cycles},
Numer. Linear Algebra Appl., 30 (2023), pp. 1--27.

\bibitem{Luo2006}
{\sc H. Luo, J. D. Baum and R. Lohner},
{\it A $p$-multigrid discontinuous Galerkin method for the Euler equations on unstructured grids},
J. Comput. Phys., 211 (2006), pp. 767--783.

\bibitem{Maday1988}
{\sc Y. Maday and R. Munoz}, 
{\it Spectral element multigrid. II. Theoretical justification},
J. Sci. Comput., 3 (1988), pp. 323--353.


\bibitem{Perugia2002}
{\sc I. Perugia and D. Schotzau}, 
{\it An $hp$-analysis of the local discontinuous Galerkin method for diffusion problems}, 
J. Sci. Comput., 17 (2002), pp. 561-571.

\bibitem{Ronquist1987}
{\sc E. M. Ronquist and A. T. Patera},
{\it Spectral element multigrid. I. Formulation and numerical results},
J. Sci. Comput., 2 (1987), pp. 389--406.

\bibitem{Stamm10}
{\sc B. Stamm and T. P. Wihler}, 
{\it hp-optimal discontinuous Galerkin methods for linear elliptic problems}, 
Math. Comp., 79 (2010), pp. 2117-2133. 

\bibitem{Elman2014}
{\sc H. C. Elman, D. J. Silvester and  A. J. Wathen}
{\it Finite elements and fast iterative solvers: with applications in incompressible fluid dynamics.}
Oxford university press, (2014).

\bibitem{Stuben1999}
{\sc K. Stuben}, 
{\it Algebraic Multigrid (AMG) : An Introduction With Applications},


\bibitem{Thiele2022}
{\sc C. Thiele and B. Riviere},
{\it $p$-multigrid with partial smoothing: An efficient
preconditioner for discontinuous Galerkin discretizations with modal bases},
J. Comput. Appl. Math., 402 (2022), pp. 113815.

\bibitem{Thiele2017}
{\sc C. Thiele, M. Araya-Polo, F. O. Alpak, B. Riviere and F. Frank},
{\it Inexact hierarchical scale separation: A two-scale approach for linear systems from discontinuous Galerkin discretizations},
Comput. Math. Appl., 74 (2017), pp. 1769--1778.

\bibitem{Vassilev2008}
{\sc P. S. Vassilevski}
{\it Multilevel block factorization preconditioners: Matrix-based analysis and algorithms for solving finite element equations},
Springer Science \& Business Media, (2008).

\bibitem{Wildey2019}
{\sc T. Wildey, S. Muralikrishnan and T. Bui-Thanh},
{\it Unified geometric multigrid algorithm for hybridized high-order finite element methods},
SIAM J. Sci. Comput., 41 (2019), pp. S172--S195.


\bibitem{Xu2024}
{\sc F. Xu, B. Wang, M. Xie},
{\it Local and parallel multigrid method for semilinear
Neumann problem with nonlinear boundary condition},
Numer. Algorithms, 96 (2024), pp. 185--210.

\bibitem{Xu2017_review}
{\sc J. Xu and L. Zikatanov},
{\it Algebraic multigrid methods},
Acta Numer., 26 (2017), pp. 591--721.

\bibitem{Xu1992}
{\sc J. Xu},
{\it Iterative methods by space decomposition and subspace correction},
SIAM review, 34 (1992), pp. 581--613.

\bibitem{Zhang2009}
{\sc L. Zhang},
{\it A Parallel Algorithm for Adaptive Local Refinement of Tetrahedral Meshes Using Bi-section},
Numer. Math. Theor. Meth. Appl., 2 (2009), pp. 65--89. 
(Version 0.9.7: http://lsec.cc.ac.cn/phg)
 
\end{thebibliography}


\newpage
\section{Appendix}
\label{app}
In this appendix, we will give the upper bound of $\rho(\bbD_k^{-1} \bbA_k)$ in Lemma \ref{Lemma_omega-k}.
Our proof proceeds in three stages, 
first, we state without proof certain properties of Legendre orthogonal polynomials.
Next, we establish the spectral radius estimate $\rho(\bbD_k^{-1} \bbA_k)$ for the one-dimensional case.
Finally, we extend these results to higher dimensions.

\subsection{ One-dimensional Legendre polynomials }
Firstly, we recall some properties of the Legendre polynomials $\{l_{\sigma} : \sigma \in \mathbb{N}\}$ on reference element $\widehat{T}=[-1,1]$.
The Legendre polynomials $\{l_{\sigma}:  \sigma \in \mathbb{N}\}$ on $\widehat{T}$ have definite parity, i.e.
\begin{align}\label{l-parity}
l_{\sigma} (-\xi)  = (-1)^{\sigma} l_{\sigma} (\xi) \quad \forall \, \xi \in \widehat{T}.
\end{align}
The parity and the normalization condition $l_{\sigma} (1) = 1,  \sigma \in \mathbb{N}$ implies $l_{\sigma} (-1)  = (-1)^{\sigma},  \sigma \in \mathbb{N}.$
The polynomials $\left\{\dot{l}_{\sigma}: \sigma \in \mathbb{N}^{+}\right\}$ 
denote the derivative of Legendre polynomials, 
those polynomials obey the Christoffel's expansion given by 
\begin{align}\label{Christoffel}
  \dot{l}_{\sigma}  = \left(2 \sigma - 1\right) l_{\sigma -1}  + \left(2 \sigma - 5\right) l_{\sigma - 3}  + \left(2 \sigma - 9\right) l_{\sigma - 5}  + \cdots,
\end{align} 
the last term of the series being $3 l_{1} $ or $l_{0}$ according as $\sigma$ is even or odd.
Combine \eqref{l-parity} and \eqref{Christoffel}, we have the definite parity for polynomials $\{\dot{l}_{\sigma} : \sigma \in \mathbb{N}^{+} \}$ as follows.
\begin{align}\label{lp-parity}
  \dot{l}_{\sigma} (-\xi)  = (-1)^{\sigma+1} \dot{l}_{\sigma} (\xi) \quad \forall \, \xi \in \widehat{T}.
\end{align}
By \eqref{Christoffel} and \eqref{lp-parity}, it is easy to calculate the pointwise evaluations as follows.
\begin{align}\label{lp-pointwise}
  \dot{l}_{\sigma} (1)  = \frac{ \sigma(\sigma+1) }{ 2 } , \quad
  \dot{l}_{\sigma} (-1)  = (-1)^{\sigma+1} \frac{ \sigma(\sigma+1) }{ 2 } .
\end{align} 
The orthogonal property of Legendre polynomials $\{l_{\sigma} , \sigma \in \mathbb{N}\}$ on $\widehat{T}$ can be written as
\begin{align}\label{l-orth}
  \SP{l_{\sigma} }{l_{\sigma'} }_{\widehat{T}} = \frac{ 2 }{ 2\sigma + 1 } \delta_{\sigma \sigma'}. 
\end{align} 
Then the Christoffel's expansion \eqref{Christoffel} can be written in compact form as 
\begin{align}\label{Christoffel-reform}
  \dot{l}_{\sigma}  = 
  \frac{ 2 }{ \NLtwo[\widehat{T}]{l_{\sigma -1} }^{2}  } l_{\sigma-1}  +
  \frac{ 2 }{ \NLtwo[\widehat{T}]{l_{\sigma -3} }^{2}  } l_{\sigma-3}  +
  \frac{ 2 }{ \NLtwo[\widehat{T}]{l_{\sigma -5} }^{2}  } l_{\sigma-5}  +
  \cdots. 
\end{align} 
Combine \eqref{l-orth} and \eqref{Christoffel-reform}, 
we can calculate the integral of $\dot{l}_\sigma, \dot{l}_{\sigma'}$  
on $\widehat{T}$ with $\sigma, \sigma' \in \mathbb{N}^+$
\begin{align}\label{lp-orth}
  \SP{\dot{l}_{\sigma} }{\dot{l}_{\sigma'} }_{\widehat{T}}  =  
  \begin{cases}
    \widetilde{\sigma} \left(\widetilde{\sigma}+1\right) \quad & \text{mod}(\sigma - \sigma',2) = 0,\\
    0 \quad & \text{mod}(\sigma - \sigma',2) = 1,
  \end{cases}
  \quad \text{where}\quad  \widetilde{\sigma}=\min(\sigma,\sigma'). 
\end{align} 

For any basis $\psi_T^\sigma \in Q_{k} (T)$ defined by \eqref{basis}, 
we have the following formulation based on affine transformation $F_{T}$ from $\widehat{T}$ to $T = [x^{l} _T,x^{r} _T]  $,
\begin{align}\label{basis1D}
\psi_{T}^{\sigma} (x) & = l_{\sigma} (\xi) \quad 
 \; \xi = F_{T}^{-1}  (x) \in \widehat{T}\quad \forall \, x \in T,
\end{align} 
where $F_{T} $ is defined by $F_{T} (\xi) : =  x^{l} _{T} + ({ h } / { 2 } ) (\xi + 1)$ and $h = x^{r} _{T} - x^{l} _{T}$. 
For any $\psi_{T} ^{\sigma} ,\psi_{T} ^{\sigma'} \in Q_{k} (T)$, 
combining the above properties \eqref{l-parity}, \eqref{lp-parity}, \eqref{lp-pointwise}, 
\eqref{l-orth}, \eqref{lp-orth} of Legendre polynomials and \eqref{basis1D}, we get
\begin{align}
  \psi_{T} ^{\sigma} (x^{r}_{T} ) & = 1,& \quad \psi_{T} ^{\sigma} (x^{l} _{T} ) & = (-1)^{\sigma} , \notag \\
  \dot{\psi_T^{\sigma}} (x^{r}_{T} ) & = \frac{ \sigma (\sigma + 1) }{ h } , & \quad 
  \dot{\psi_T^{\sigma}} (x^{l} _{T} ) &= (-1)^{\sigma+1} \frac{ \sigma (\sigma + 1) }{ h } , \label{props} \\
  \SP{\psi_{T}^{\sigma}}{\psi_{T} ^{\sigma'} }_{T} & = \frac{ h }{ 2\sigma+1 } \delta_{\sigma\sigma'}, & \quad
\SP{\dot{\psi_T^{\sigma}}}{\dot{\psi_{T} ^{\sigma'}}}_{T} & =
\begin{cases}
   2\widetilde{\sigma} \left(\widetilde{\sigma} + 1\right) / h \quad 
  & \text{mod} (\sigma - \sigma', 2) = 0, \\
  0 \quad 
  & \text{mod}(\sigma - \sigma',2) = 1.
\end{cases} \notag
\end{align}

\subsection{ Spectral estimation on the one-dimensional case}
Next, we start from the one-dimensional case with broken polynomial spaces 
$\{U_k(\Ct_h)\}_{k=1}^{p}  $ defined by \eqref{DGspace}. 
For a given element $T = [x_{T} ^{l} , x_{T} ^{r} ]  \in \Ct_{h}$, we define the elementary stencil $\Cs \left(\Ca_{k} , T; \Ct_{h} \right)$ associated with the bilinear form $\Ca_{k}$ as
\begin{align*} 
\Cs \left(\Ca_{k} , T ; \Ct_{h} \right) := 
\left\{
T' \in \Ct_{h}: \;
\exists \psi_{T} ^{\sigma} \in Q_{k} (T), \;
\exists \psi_{T'} ^{\sigma'} \in Q_{k} (T'), \;
\Ca_{k} (\psi_{T} ^{\sigma} , \psi_{T'} ^{\sigma'} ) \neq 0
\right\}.
\end{align*}
By the definition \eqref{IPDG-a} of $\Ca_{k} \left(\cdot, \cdot \right)$ in the one-dimensional case, 
we get $ \text{card} \left\{ \Cs \left( \Ca_{k} , T; \Ct_{h} \right) \right\} \le 3$,
where $ \text{card} \left\{ V \right\}$ denotes the total number of elements in $V$.
Suppose $T' = T = [x_T^l, x_T^r]$, then we have 
\begin{equation*}
\begin{aligned}
\Ca_{k}(\psi_T^\sigma, \psi_{T}^{\sigma'}) 
&=& &\int_{T} \partial_x \psi_T^{\sigma} \partial_x \psi_T^{\sigma'} dx 
- \left( \avg{\partial_x \psi_T^{\sigma}(x_T^r) } \jump{\psi_T^{\sigma'}(x_T^r) }  
+ \avg{ \partial_x \psi_T^{\sigma'} (x_T^r)}  \jump{\psi_T^{\sigma} (x_T^r) } \right)\\
& & &-\left( \avg{\partial_x \psi_T^{\sigma} (x_T^l) } \jump{ \psi_T^{\sigma'} (x_T^l) }
+ \avg{\partial_x \psi_T^{\sigma'} (x_T^l ) } \jump{ \psi_T^{\sigma} (x_T^l) }  \right) \\
& & &+\frac{\alpha_0 k^2}{h} \left( \jump{ \psi_T^{\sigma} (x_T^r) } \jump{ \psi_T^{\sigma'} (x_T^r)} 
+  \jump{ \psi_T^{\sigma} (x_T^l) } \jump{\psi_T^{\sigma'} (x_T^l)} \right).
\end{aligned}
\end{equation*}
Since the basis function $\psi_T^\sigma$ vanishes on all elements except $T$, 
we conclude that,
\begin{align*}
&\jump{ \psi_T^{\sigma} (x_T^l) } = -\psi_T^{\sigma} (x_T^l) = (-1)^{\sigma+1}, \quad
&&\jump{ \psi_T^{\sigma} (x_T^r) } = \psi_T^{\sigma} (x_T^r) = 1, \\
&\avg{\partial_x \psi_T^{\sigma}(x_T^l) } = \frac12 \partial_x \psi_T^{\sigma}(x_T^l) = (-1)^{\sigma+1} \frac{\sigma(\sigma+1)}{2h}, \quad
&&\avg{\partial_x \psi_T^{\sigma}(x_T^r) } = \frac12 \partial_x \psi_T^{\sigma}(x_T^r) = \frac{\sigma(\sigma+1)}{2h} .
\end{align*}
Applying these properties into the above bilinear form. 
When $\text{mod} \left(\sigma - \sigma', 2\right) = 1$, due to parity mismatch, 
$\Ca_{k} (\psi_{T}^{\sigma} , \psi_{T} ^{\sigma'} ) = 0$ vanishes.
When $\text{mod} \left(\sigma - \sigma', 2\right) = 0$, we have 
\begin{align*}
\Ca_{k} (\psi_{T}^{\sigma} , \psi_{T} ^{\sigma'} )
&= \frac{2\tilde{\sigma}(\tilde{\sigma} + 1)}{h} - \left(1 + (-1)^{\sigma+\sigma'} \right) \frac{\sigma(\sigma+1) + \sigma'(\sigma'+1)}{2h} 
+ \frac{\alpha_0 k^2}{h} \left(1 + (-1)^{\sigma+\sigma'} \right) \\
&= \frac{1}{ h } \left[ 
    2 \alpha_0 k^{2} 
    +  \widetilde{\sigma}(\widetilde{\sigma} + 1) 
    -  \bar{\sigma}(\bar{\sigma} + 1)  
  \right] 
\end{align*}
where $\bar{\sigma} = \max(\sigma, \sigma'), \tilde{\sigma} = \min(\sigma, \sigma')$.
For $T' \neq T$, suppose $T'$ is the right adjacent cell of $T$ with $T\bigcap T' = x_T^r$.
Then we have,
\begin{align*}
\Ca_{k}(\psi_T^\sigma, \psi_{T'}^{\sigma'}) 
=& - \left( \avg{\partial_x \psi_T^{\sigma}(x_T^r) } \jump{\psi_{T'}^{\sigma'}(x_T^r) }   
+ \avg{ \partial_x \psi_{T'}^{\sigma'} (x_T^r)}  \jump{\psi_T^{\sigma} (x_T^r) } \right) \\
 & +\frac{\alpha_0 k^2}{h} \left( \jump{ \psi_T^{\sigma} (x_T^r) } \jump{ \psi_{T'}^{\sigma'} (x_T^r)} \right) \\
=&(-1)^{\sigma' + 1}\left(\frac{\alpha_0 k^2}{h} - \frac{\sigma(\sigma + 1)}{2h} - \frac{\sigma' (\sigma' + 1)}{2h} \right).
\end{align*}

Therefore, we come to the conclusion that, for any $T' \in \Cs (\Ca_{k} , T; \Ct_{h} )$, we have
\begin{align}
\Ca_{k} (\psi_{T}^{\sigma} , \psi_{T} ^{\sigma'} ) =   
& \begin{cases}
  \dfrac{1}{ h } \left[ 
    2 \alpha_0 k^{2} 
    +  \widetilde{\sigma}(\widetilde{\sigma} + 1) 
    -  \bar{\sigma}(\bar{\sigma} + 1)  
  \right] \quad
  & \text{mod} \left(\sigma - \sigma', 2\right) = 0, \\
  0 \quad & \text{mod} \left(\sigma - \sigma', 2\right) = 1,
\end{cases} \label{att}\\
\SN{\Ca_{k} (\psi_{T}^{\sigma} , \psi_{T^{'} } ^{\sigma'} )} =  
& \frac{1}{ 2 h } 
  \left| 
    2 \alpha_{0} k^{2}
    -  {\sigma}({\sigma} + 1) 
    -  {\sigma'}({\sigma'} + 1) 
  \right| \quad 
  \text{if} \quad T' \neq T. \label{attp} 
\end{align}

For a given basis $\psi_{T} ^{\sigma} \in Q_{k} (T)$, we introduce the set $\Cn(\Ca_{k} , \psi_{T} ^{\sigma} ; Q_{k} (T')) \subset \mathbb{N}$ associated with the bilinear form $\Ca_{k} $, given by
\begin{align*}
  \Cn(\Ca_{k} , \psi_{T} ^{\sigma} ; Q_{k} (T')) :=
  \left\{ 
  \sigma' \in \mathbb{N}:
  \Ca_{k} (\psi_{T} ^{\sigma} , \psi_{T'} ^{\sigma'} ) \neq 0,
  \forall \, \psi_{T'} ^{\sigma'} \in Q_{k} (T')
  \right\} \quad 
  \forall \, T' \in \Cs (\Ca_{k} , T; \Ct_{h} ).
\end{align*}
By \eqref{att} and \eqref{attp}, we find that
\begin{align}\label{NCn}
  \text{card} \left\{ \Cn (\Ca_{k},\psi_{T} ^{\sigma} ; Q_{k} (T')) \right\} \le
  k + 1 \quad
  \forall \, T' \in \Cs (\Ca_{k} , T; \Ct_{h} ).
\end{align}
For any $u_k\in U_k(\Ct_h)$, 
the vector of unknowns is denoted by $\Vu_{k} $, 
consisting of coefficients of the expansion 
$u_k = \sum_{T\in \Ct_h}\sum_{\sigma=0}^k u_T^{\sigma}  \psi_T^{\sigma} $.

Now, we prove the upper bound of $\rho \left(\bbD_{k} ^{-1} \bbA_{k} \right) $ in the one-dimensional case.
By the coercivity and continuity of $\Ca_{k} (\cdot, \cdot)$ in Lemma \ref{lem2.1}, we have
\begin{align}\label{ccle}
\begin{split}
\Vu_{k} ^\top \mathbb{A}_k \Vu_{k} 
= & 
\Ca_k \Big(
  \sum_{T \in \Ct_{h} } \sum_{\SN{\sigma} \le k} u_{T} ^{\sigma}  \psi_{T} ^{\sigma} , 
  \sum_{T' \in \Ct_{h} }  \sum_{\SN{\sigma'} \le k} u_{T'}^{\sigma'} \psi_{T'}^{\sigma'} 
  \Big) \\
\le & 
  \sum_{T \in \Ct_{h} } \sum_{\SN{\sigma} \le k} \SN{u_{T} ^{\sigma}}  
  \sum_{T' \in \Ct_{h} }  \sum_{\SN{\sigma'} \le k} \SN{u_{T'}^{\sigma'} }
  \SN{\Ca_{k} ( \psi_{T}^{\sigma} , \psi_{T'}^{\sigma'} )}
\\ 
\lesssim &
  \sum_{T \in \Ct_{h} } \sum_{\SN{\sigma} \le k} \SN{u_{T} ^{\sigma}}  
  \Ca^{\frac12}_{k} ( \psi_{T}^{\sigma} , \psi_{T}^{\sigma} )
  \sum_{T' \in \Ct_{h} }  \sum_{\SN{\sigma'} \le k} \SN{u_{T'}^{\sigma'} }
  \Ca^{\frac12}_{k} ( \psi_{T'}^{\sigma'} , \psi_{T'}^{\sigma'} ) 
\\ 
\lesssim & 
\sum_{T \in \Ct_{h} } \sum_{\SN{\sigma} \le k} 
\SN{u_{T} ^{\sigma} } 
\Ca^{\frac{1}{2}}_k(\psi_{T} ^{\sigma} , \psi_{T}^{\sigma} )
\sum_{T' \in \Cs( \Ca_{k} , T; \Ct_{h} )} 
\sum _{\sigma' \in \Cn(\Ca_{k} , \psi_{T} ^{\sigma} ; Q_{k} (T'))} 
\SN{u_{T'}^{\sigma'} }
\Ca^{\frac{1}{2}}_k(\psi_{T'} ^{\sigma'} , \psi_{T'}^{\sigma'} )
\end{split}
\end{align}
where $\Ca_k(\psi_{T'} ^{\sigma'} , \psi_{T'}^{\sigma'} )\neq 0$, 
only when $T' \in \Cs( \Ca_{k} , T; \Ct_{h} )$ and $\sigma' \in \Cn(\Ca_{k} , \psi_{T} ^{\sigma} ; Q_{k} (T'))$. 
Notice that $\text{card} \left\{ \Cs \left( \Ca_{k} , T; \Ct_{h} \right) \right\} \le 3$ and 
$\text{card} \left\{ \Cn (\Ca_{k},\psi_{T} ^{\sigma} ; Q_{k} (T')) \right\} \le k + 1$, 
for given $T \in \Ct_{h}$ and $\sigma \le k$, using Schwarz's inequality, we have 
\begin{align}\label{si-1}
\begin{split}
& \sum_{T' \in \Cs( \Ca_{k} , T; \Ct_{h} )} 
\sum _{\sigma' \in \Cn(\Ca_{k} , \psi_{T} ^{\sigma} ; Q_{k} (T'))} 
\SN{u_{T'}^{\sigma'} }
\Ca^{\frac{1}{2}}_k(\psi_{T'} ^{\sigma'} , \psi_{T'}^{\sigma'} ) \\ 
\leq & 
\textcolor{blue}{
\Big(\sum_{T' \in \Cs( \Ca_{k} , T; \Ct_{h} )} 
\sum _{\sigma' \in \Cn(\Ca_{k} , \psi_{T} ^{\sigma} ; Q_{k} (T'))}  
1^2 \Big)^{\frac{1}{2}}  \times
}  \\
& \Big(\sum_{T' \in \Cs( \Ca_{k} , T; \Ct_{h} )} 
\sum _{\sigma' \in \Cn(\Ca_{k} , \psi_{T} ^{\sigma} ; Q_{k} (T'))}  
|u_{T'}^{\sigma'}|^2 \Ca_k(\psi_{T'} ^{\sigma'} , \psi_{T'}^{\sigma'} ) 
\Big)^{\frac{1}{2}}  \\
\leq &
\sqrt{3(k+1)} 
\Big(\sum_{T' \in \Cs( \Ca_{k} , T; \Ct_{h} )} 
\sum _{\sigma' \in \Cn(\Ca_{k} , \psi_{T} ^{\sigma} ; Q_{k} (T'))}  
|u_{T'}^{\sigma'}|^2 \Ca_k(\psi_{T'} ^{\sigma'} , \psi_{T'}^{\sigma'} ) 
\Big)^{\frac{1}{2}}.
\end{split}
\end{align}
It should be noted that $\sum_{T\in \Ct_h} \sum_{\SN{\sigma} \leq k } 
\SN{u_{T} ^{\sigma} }^2 \Ca_k(\psi_{T} ^{\sigma} , \psi_{T}^{\sigma} ) = \Vu_{k} ^{\top} \bbD_{k} \Vu_{k} $,
then using \eqref{ccle} and \eqref{si-1}, we have
\begin{align}\label{si-2}
\begin{split}
\Vu_{k}^{\top} \bbA_{k} \Vu_{k} \lesssim 
& \sqrt{3(k+1)} 
\sum_{T \in \Ct_{h} } \sum_{\SN{\sigma} \le k} 
\SN{u_{T} ^{\sigma} } 
\Ca^{\frac{1}{2}}_k(\psi_{T} ^{\sigma} , \psi_{T}^{\sigma} ) \times \\
& \Big(
\sum_{T' \in \Cs( \Ca_{k} , T; \Ct_{h} )} 
\sum _{\sigma' \in \Cn(\Ca_{k} , \psi_{T} ^{\sigma} ; Q_{k} (T'))}  
|u_{T'}^{\sigma'}|^2 \Ca_k(\psi_{T'} ^{\sigma'} , \psi_{T'}^{\sigma'} ) 
\Big)^{\frac{1}{2}}, \\
\lesssim 
& \sqrt{3(k+1)} \left( 
\Vu_{k} ^{\top} \bbD_{k} \Vu_{k} 
\right)^{\frac{1}{2}} \times \\
& \Big( \sum_{T \in \Ct_{h} } \sum_{\SN{\sigma} \le k} 
\sum_{T' \in \Cs( \Ca_{k} , T; \Ct_{h} )} 
\sum _{\sigma' \in \Cn(\Ca_{k} , \psi_{T} ^{\sigma} ; Q_{k} (T'))}  
|u_{T'}^{\sigma'}|^2 \Ca_k(\psi_{T'} ^{\sigma'} , \psi_{T'}^{\sigma'} ) 
 \Big)^{\frac{1}{2}},
\end{split}
\end{align}
where the last inequality is obtained from the Schwarz's inequality.
By the symmetry of the bilinear form, 
$T' \in \Cs( \Ca_{k} , T; \Ct_{h} ), \sigma' \in \Cn(\Ca_{k} , \psi_{T} ^{\sigma} ; Q_{k} (T'))$ yields
$T \in \Cs( \Ca_{k} , T'; \Ct_{h} ), \sigma \in \Cn(\Ca_{k} , \psi_{T'} ^{\sigma'} ; Q_{k} (T))$.
Swapping the sum sign order and using the property \eqref{NCn}, we have
\begin{align}\label{swapsum}
\begin{split}
&\sum_{T \in \Ct_{h} } \sum_{\SN{\sigma} \le k} 
\sum_{T' \in \Cs( \Ca_{k} , T; \Ct_{h} )} 
\sum _{\sigma' \in \Cn(\Ca_{k} , \psi_{T} ^{\sigma} ; Q_{k} (T'))}  
|u_{T'}^{\sigma'}|^2 \Ca_k(\psi_{T'} ^{\sigma'} , \psi_{T'}^{\sigma'} ) \\
= &\sum_{T' \in \Ct_{h} } \sum_{\SN{\sigma'} \le k} 
\sum_{T \in \Cs( \Ca_{k} , T'; \Ct_{h} )} 
\sum _{\sigma \in \Cn(\Ca_{k} , \psi_{T'} ^{\sigma'} ; Q_{k} (T))}  
|u_{T'}^{\sigma'}|^2 \Ca_k(\psi_{T'} ^{\sigma'} , \psi_{T'}^{\sigma'} ) 
\leq 3(k+1)  \Vu_{k} ^{\top} \bbD_{k} \Vu_{k}. 
\end{split}
\end{align}
Using \eqref{si-2} and \eqref{swapsum}, 
we get the estimation ${\Vu_{k} ^\top \bbA_k \Vu_{k}} / {\Vu_{k} ^\top \bbD_k \Vu_{k}  }\lesssim k + 1$.
Then
\begin{align*}
\rho(\bbD^{-1}_k \bbA_k) 
= & \rho(\bbD^{\frac{1}{ 2 }} \bbD_{k} ^{-1} \bbA_{k} \bbD_{k} ^{- \frac{1}{ 2 } } )
= \rho(\bbD^{-\frac{1}{2}}_k \bbA_k \bbD^{-\frac{1}{2}}_k) \\
=  & \sup\limits_{\Vv_{k} \neq \boldsymbol{0}}
\frac{ \Vv_{k} ^\top \bbD_{k} ^{-\frac{1}{2}} \bbA_k \bbD^{-\frac{1}{2}}_k\Vv_{k}  }{ \Vv_{k} ^\top \Vv_{k}  } \\
= & \sup \limits_{\substack{\Vu_{k} = \bbD_{k} ^{-\frac{1}{ 2 } } \Vv_{k} \\ \Vv_{k}  \neq 0  }} 
\frac{\Vu_{k} ^\top \bbA_k \Vu_{k}  }{\Vu_{k} ^\top \bbD_k \Vu_{k}  }
\lesssim k + 1.
\end{align*}

\begin{remark}
The core of our spectral radius estimate lies in the key observation that the number of non-zero entries scales as $\Co(k)$ 
rather than $\Co(k^d)$. 
While this is intuitively reasonable in 1D settings with \eqref{NCn}, 
we will prove its validity in higher dimensions, which is a critical improvement 
that enhances the smoothing properties of our $p$-multigrid method.
\end{remark}

\subsection{Extension to higher-dimensional case}
For $d$-dimensional case, we have the following tensor-product formulation of Legendre polynomial in $Q_{k} (T)$, where the element $T$ is given by $T = \Pi_{i=1} ^{d}  [x_{i} ^{l} , x_{i} ^{r} ] $,
\begin{align*}
\Psi_{T} ^{\sigmabf} \left(\Bx\right) = \prod_{i=1}^{d} \psi_{T} ^{\sigma_{i} } \left(x_{i} \right) 
\quad \text{with} \quad \psi_{T} ^{\sigma_{i} } (x_{i} ) = l_{\sigma_{i} } (\xi_{i} ) \quad 
\; \xibf = \BF_{T} ^{-1} \left(\Bx\right) \quad \forall \, \Bx \in T.
\end{align*}
Here $\BF_{T} \left(\xibf\right) = \left(F_{T} ^{i} (\xi_{i} ) : i = 1, \dots , d\right)^{\top}  $ is the affine transformation from $\widehat{T} = [-1,1]^{d} $ to $T$,
and $F_{T} ^{i} (\xi_{i} ) := x_{i}^{l} +  \left(h / 2\right) (\xi_{i} + 1)$.

For any $T' \in \Cs (\Ca_{k} , T; \Ct_{h} ) $,
suppose $T'\neq T$ and $T\bigcap T' = F$ is a plane with $x = x_i^l$, then the bilinear form is 
\begin{align*}
\Ca_{k} (\Psi_{T} ^{\sigmabf}, \Psi_{T'} ^{\sigmabf'}) 
:=& - \DP{\avg{ \nabla \Psi_{T} ^{\sigmabf} }\cdot\Bn }{\jump{ \Psi_{T'} ^{\sigmabf'} }}_{F} 
- \DP{\avg{\nabla \Psi_{T'} ^{\sigmabf'} }\cdot \Bn }{\jump{ \Psi_{T} ^{\sigmabf} } }_{F} \\
&+ \frac{\alpha_0 k^2}{h} \DP{\jump{ \Psi_{T} ^{\sigmabf} }}{\jump{ \Psi_{T'} ^{\sigmabf'} }}_{F}.
\end{align*}
The computation along the $x$-direction follows analogous procedures to the 1D case previously established. 
For conciseness, we omit the technical derivations and state the final result:
\begin{align*}
  \Ca_{k} (\Psi_{T} ^{\sigmabf} , \Psi_{T} ^{\sigmabf'}) = &
  \begin{cases}
  h^{d - 2} 
  \sum\limits_{i = 1} ^{d} 
  \left[
    2 \alpha_0 k^{2} 
    + \widetilde{\sigma}_{i} (\widetilde{\sigma}_{i} + 1) 
    - \bar{\sigma}_{i} (\bar{\sigma}_{i}  + 1)
  \right]
  \prod\limits_{\substack{ j = 1 \\ j \neq i}}^{d} \dfrac{ \delta_{\sigma_{j} \sigma'_{j} }  }{ 2 \sigma_{j} + 1 } \quad
  \text{mod}(\sigma_{i} - \sigma'_{i} , 2 ) = 0,\\ 
  0 \hfill \exists \; i \in \{1, 2, \dots, d\}, \text{mod}( \sigma_{i} - \sigma'_{i} , 2 ) = 1,  
  \end{cases} \\
  \Ca_{k} (\Psi_{T} ^{\sigmabf} , \Psi_{T'} ^{\sigmabf'}) = &
  \frac{ 1 }{ 2h } 
  \left[
    2 \alpha_{0} k^{2} 
    - \sigma_{i} (\sigma_{i} + 1)
    - \sigma'_{i} (\sigma'_{i} + 1)
  \right]
  \prod_{\substack{j = 1\\ j \neq i}}^{d}  \frac{ \delta_{\sigma_{j} \sigma'_{j} }  }{ 2 \sigma_{j} + 1 } \quad
  \text{for} \quad T' \neq T . 
\end{align*}
By the above identities, 
we have the additional properties 
\begin{align}
  \text{card} \left\{ \Cs (\Ca_{k},{T} ;\Ct_{h}  )\right\} & \le 
  2d+1 \quad 
  \;\;\;\forall \, T \in \Ct_{h},  \label{nS} \\
  \text{card} \left\{ \Cn (\Ca_{k},\Psi_{T} ^{\sigmabf} ; Q_{k} (T')) \right\} & \le
  k + 1 \quad
  \quad \; \, \forall \, T' \in \Cs (\Ca_{k} , T; \Ct_{h} ). \label{nN}
\end{align}
By using \eqref{nS}, and \eqref{nN},
the one-dimensional proof of the spectral radius can be extended to $d$-dimensional cases.
We can also get $\rho(\bbD^{-1}_k\bbA_k) \lesssim k + 1$. 
The details are omitted here.


\end{document}